\documentclass[hidelinks,onefignum,onetabnum]{siamart251216}

\usepackage{lipsum}
\usepackage{amsfonts}
\usepackage{graphicx}
\usepackage{epstopdf}
\usepackage{algpseudocode}
\usepackage{amsmath}
\usepackage{mathrsfs}
\usepackage{cases}
\usepackage{subcaption}
\usepackage{siunitx}
\usepackage{booktabs}

\ifpdf
  \DeclareGraphicsExtensions{.eps,.pdf,.png,.jpg}
\else
  \DeclareGraphicsExtensions{.eps}
\fi

\newsiamremark{remark}{Remark}
\newsiamremark{hypothesis}{Hypothesis}
\crefname{hypothesis}{Hypothesis}{Hypotheses}
\newsiamthm{claim}{Claim}
\newsiamremark{fact}{Fact}
\crefname{fact}{Fact}{Facts}

\headers{The scattering law for a RTKM}{S. Lei, J. Liang, W. Sun and M. Tang}

\title{The scattering law at dark-light interfaces\\ for a run-and-tumble kinetic model\thanks{Submitted to the editors DATE.
\funding{This work was funded by the Fog Research Institute under contract no.~FRI-454.}}}

\author{Shicong Lei \thanks{School of Mathematical Sciences, Shanghai Jiao Tong University, Shanghai 200240, P. R. China (\email{leishicong322@sjtu.edu.cn})}
\and Jiangyan Liang \thanks{School of Mathematics and Statistics, Wuhan University of Technology, Wuhan 430070, P. R. China (\email{jyliang@whut.edu.cn}).}
\and Weiran Sun \thanks{Department of Mathematics, Simon Fraser University, 8888 University Dr., Burnaby, BC V5A 1S6, Canada (\email{weirans@sfu.ca}).}
\and Min Tang \thanks{School of Mathematics and Institute of Natural Sciences, MOE-LSC, Shanghai Jiao Tong University, Shanghai 200240, P. R. China (\email{tangmin@sjtu.edu.cn}).}
}

\usepackage{amsopn}

\usepackage{amssymb}

\def\p{\partial}

\def\d{\mathrm{d}}
\def\no{\nonumber}
\def\R{\mathbb{R}}

\def\eps{\varepsilon}

\def\'{\prime}

\newcommand{\vpran}[1]{\left(#1\right)}
\newcommand{\nn}{\nonumber}

\ifpdf
\hypersetup{
  pdftitle={The scattering law at dark-light interfaces for a run-and-tumble kinetic model},
  pdfauthor={Shicong Lei, Jiangyan Liang, Weiran Sun and Min Tang}
}
\fi

\title{The interface scattering law for run-and-tumble kinetic models \thanks{\funding{The author J. Liang is supported by NSFC 12601959. W. Sun is partially supported by NSERC R832717. M. Tang is funded by NSFC U2530207.}}}

\author{Shicong Lei \thanks{School of Mathematical Sciences, Shanghai Jiao Tong University, Shanghai 200240, P. R. China (\email{leishicong322@sjtu.edu.cn})}
\and Jiangyan Liang \thanks{School of Mathematics and Statistics, Wuhan University of Technology, Wuhan 430070, P. R. China (\email{jyliang@whut.edu.cn}).}
\and Weiran Sun \thanks{Department of Mathematics, Simon Fraser University, 8888 University Dr., Burnaby, BC V5A 1S6, Canada (\email{weirans@sfu.ca}).}
\and Min Tang \thanks{School of Mathematics and Institute of Natural Sciences, MOE-LSC, Shanghai Jiao Tong University, Shanghai 200240, P. R. China (\email{tangmin@sjtu.edu.cn}).}
}

\headers{Interface scattering law for RTKM}{S. Lei, J. Liang, W. Sun and M. Tang} 

\begin{document}

\maketitle

\begin{abstract}
For microorganisms with run-and-tumble motility, kinetic models provide a well-established description of transport in smoothly varying environments. Sharp signal changes, however, can induce asymmetric interface responses.  A motivating example is the response to dark–light boundaries:  cells moving into the light largely maintain their direction, whereas those entering the dark are effectively scattered back. Through an interface-layer analysis of pathway-based kinetic equations, we formally derive effective interface conditions for the limiting run-and-tumble model. These velocity-nonlocal conditions capture non-adiabatic scattering near signal discontinuities, linking intracellular sensory dynamics to transport across sharp environmental boundaries. For a two-velocity model, they reduce to an explicit transmission matrix. We assess the interface law numerically using biologically motivated parameters for E. coli chemotaxis.
\end{abstract}



\section{Introduction}\label{Sec:Intro}
Some bacteria or algae move by alternating runs and tumbles, a locomotion mode for navigating environmental cues such as chemicals or light.
A run is a sustained,  directional movement. A tumble is a brief reorientation, which ends when a new run starts in a random direction. This cycle lets bacteria sample their surroundings, and it plays a fundamental role in enabling effective navigation in heterogeneous environments \cite{Alvarez2014, Berg1972, Block1982, Shimizu2010, Zaferani2023}.

Kinetic transport equations have emerged as a well-accepted framework for modeling bacterial run-and-tumble behavior \cite{Wolfgang-JMB-1980,Hillen-SIAM-2000,Othmer-JMB-1988}. In this class of models, all cells are assumed to sense and respond to environmental cues uniformly. The signal tunes the tumbling frequency: cells reduce tumbling when moving in a favorable direction (e.g., toward higher nutrient concentrations or optimal light), while increasing tumbling in adverse directions \cite{Chalub-2004,Perthame-2007}. However, run-and-tumble kinetic models (RTKMs) consider scenarios where environmental signals change slowly and smoothly.

It is now well established that fast environmental changes demand accounting for cellular heterogeneity—different cells may respond differently to the same external stimulus \cite{Bray-1998,Waite-2018}. To capture this variability, researchers turn to pathway-based kinetic models (PBKM), which incorporate an additional internal state variable to characterize cell-specific traits \cite{Erban-2004,Erban-2005}. For example, in {\it E. coli}, this variable could represent the methylation level of chemoreceptors \cite{Jiang2010SPECS,Si2014}; in other microbes, it might correspond to the expression level of light-sensitive proteins. This internal variable tracks the cell's dynamic adaptation to signals, enabling the model to capture diverse response patterns across a population and accurately describe behavior in rapidly changing, heterogeneous environments \cite{Si2012,Si2014}. These two classes of models were primarily proposed for {\it E. coli} chemotaxis, while other microorganisms exhibiting similar locomotion patterns can be modeled in a similar way.

Most existing experimental studies and modeling efforts focus primarily on scenarios where environmental signals change continuously \cite{Frangipane-eLife-2018,Jékely-2009,Kantsler-PNAS-2013}, while cells exhibit distinct dynamics when encountering abrupt signal discontinuities. Recent experimental observations have revealed that certain phototactic bacteria or algae exhibit striking differential responses when crossing dark-light interfaces: cells moving from dark to light regions (i.e., in the direction of increasing light intensity) retain their original direction of movement largely unchanged, sustaining their run phase to exploit the favorable light environment. In contrast, cells transitioning from light to dark domains (moving toward decreasing light) are effectively scattered back from the dark—their movement is redirected away from the interface, with a strong bias toward re-entering the light region \cite{Polin-PRL-2015,Drescher-PRL-2009}. To quantify this behavior, several experiments have been conducted to measure the ``scattering law" at such interfaces \cite{Leonardo-PNAS-2025}. Specifically, these studies track cells moving with an initial direction $\mathbf{v}$ and measure the probability that the cell, upon encountering the interface, scatters into a new direction $\mathbf{v}'$, thus providing the ``scattering law" at the interface \cite{Buchner-PNAS-2021,Knut-PNAS-2010}. 

At these abrupt interfaces, the RTKM fails to capture the observed directional bias and scattering dynamics. In contrast, the PBKM remains valid in this context: by integrating intracellular signaling pathway dynamics, they can replicate how cells process sudden signal shifts in real time. This raises a natural question: Can we derive the interface scattering law of the run-and-tumble kinetic transport model from the more detailed pathway-based kinetic model? Such a derivation would bridge the gap between microscopic intracellular signaling and macroscopic scattering rules \cite{Perthame-Tang-Vauchelet,Sun-Tang-2017}.

The PBKM takes the following form
 \begin{align}\label{sys:kinetic-transport model}\tag{PBKM}
    \begin{cases}
  \p_t p + v \cdot \nabla_x p + \p_m (f(m-m_c)p) = \Lambda (m-m_c)(\tfrac{1}{|V|}\int_V p \,\d v' - p),\\[2mm]
  p (0,x,v,m) = p^{in}(x,v,m),\\[2mm]
      p(t,x,v,M_l) = p(t,x,v,M_r) = 0,
\end{cases}
\end{align}
where $p(t,x,v,m)$ denotes the probability density function of bacteria at time $t$, position $x\in \mathbb{R}^d$, moving with direction $v\in V$ ($V \subseteq \R^d$ is bounded),
and with internal state $m\in [M_l,M_r]$.  The spatially dependent function $m_c(x)$ represents the equilibrium intracellular state level. It is a function of the extracellular chemical concentration (in the case of chemotaxis) or light intensity (in the case of phototaxis).
The function $f(\cdot)$ models intracellular chemical reactions and satisfies the following properties:
\begin{align}\label{eq:f}
  f(0)=0,\qquad f(y)<0 \quad\text{for } y>0 \quad\text{and}\quad f(y)>0 \quad\text{for } y<0.
\end{align}
The tumbling frequency $\Lambda$ is a monotonically increasing function of the deviation $m-m_c(x)$ of the internal state $m$ from its equilibrium state $m_c(x)$. 

Meanwhile, the RTKM ignores the internal variable and models the motion of bacteria at the cellular scale. It governs the marginal density $\bar{p} = \int_{[M_l, M_r]} p(t, \mathbf{x},\mathbf{v},m)\,\mathrm{d}m$ in the form 
\begin{equation*}\tag{RTKM}
  \partial_t \bar{p} + v \cdot\nabla_{\mathbf{x}}\bar{p}
  = \mathcal{T}[D_t m_c](\bar{p}),
  \label{eq:kinetic}
\end{equation*}
with the tumbling operator
    \begin{align*}
         \mathcal{T}[D_t m_c](\bar{p})
  = &\int_V \Bigl[
    T(D_t m_c(t,x,v'),v,v')\,\bar{p}(t, x,v') \no\\
    &- T(D_t m_c(t,x,v),v',v)\,\bar{p}(t, x,v)
  \Bigr]\,\mathrm{d}v'.
    \end{align*}
Here the path-wise gradient of the signal is defined as
\begin{equation*}
  D_t m_c(t,x,v) = \partial_t m_c(t,x) + v\cdot\nabla_{x} m_c(t,x).
  \label{eq:DtM}
\end{equation*}
Thus, in this reduced model, the tumbling kernel $T$ depends on the path-wise gradient 
instead of the signal itself.

The main difference between the PBKM and the RTKM lies in how environmental signals enter the tumbling response. In the PBKM, the response depends on the mismatch between the internal state $m$ and the local equilibrium level $m_c$, whereas in the RTKM it depends on the signal gradient $\nabla_x m_c(x)$ for a stationary $m_c$. In \cite{Perthame-Tang-Vauchelet}, the authors rigorously derived the RTKM as an asymptotic limit of the PBKM under fast adaptation, modeled by $\frac{1}{\varepsilon}f(m-m_c)$, and a stiff tumbling response, modeled by $\Lambda((m-m_c)/\varepsilon)$. As $\varepsilon\to0$, the internal-state distribution concentrates into a Dirac mass, and the signal gradient enters the limiting tumbling kernel. This derivation assumes that $m_c(x)$ varies smoothly in space. 

When $m_c$ is discontinuous, the analysis no longer applies directly near the interface. Two features become relevant. First, differences in cellular responses to abrupt signal changes suggest that the internal-state distribution in $m$ cannot be neglected \cite{Korobkova-Nature-2004,Waite-2018}. Second, cells crossing from dark to light environments can exhibit strongly reduced tumbling, whereas those crossing in the opposite direction exhibit increased tumbling. This directional asymmetry motivates assigning different asymptotic orders to $\Lambda(m-m_c)$ for different internal states \cite{Block1982,Bernardo-PNAS-2005,Victor-PNAS-2022}.
The central challenge is therefore to identify a physically meaningful scaling suitable for interface-layer analysis. To make these scaling choices explicit, we discretize the internal variable $m$ into a finite set of states. The resulting finite-state PBKM allows us to assign different asymptotic orders to the tumbling rates associated with distinct internal states near the interface.

In this work, we consider two regions with distinct constant equilibrium methylation levels separated by a sharp interface, so that the RTKM can be used to describe the interior dynamics in each region. Through an interface-layer analysis of the discrete PBKM, we derive interface conditions coupling the two RTKMs. These conditions take an integral form in velocity space. For a simplified two-flux model, we obtain an explicit interface transmission matrix.

This paper is organized as follows. Section \ref{Sec2:model} introduces the discrete model and the interface problem. Section \ref{sec:Main results} presents the main asymptotic analysis to derive the limiting system together with the interface conditions. Section \ref{sec:numerical-validation} presents numerical simulations that validate the theoretical findings and demonstrate the model's applicability to a continuous internal variable. Finally, we conclude in Section \ref{sec:conclusion}.

\section{The discrete model}\label{Sec2:model}

\subsection{The multiple species kinetic model}\label{subsec2:the kinetic model with multiple discrete species}

To simplify the notations, we consider the PBKM in slab geometry such that
\begin{equation}\label{PBKM1d}
  \p_t p +v \p_x p + \p_m (f(m-m_c)p) = \Lambda (m-m_c)\big(\frac{1}{2}\int_{-1}^1 p \,\d v' - p\big),
\end{equation}
where $v\in[-1,1]$ now represents the cosine of the angle between the particle’s movement direction and the positive $x$-axis.
The physical domain is composed of a left domain $\Omega_l$ whose equilibrium methylation level is $m_c(x) = M_l$, and a right domain $\Omega_r$ whose equilibrium methylation level is $m_c(x) = M_r$. Here $M_l< M_r$ are two different constants.

In this setting, where the methylation level $M_l\leq m\leq M_r$, we introduce the discretized version of \eqref{PBKM1d} by discretizing the continuous internal state $m$ into a finite set of $n+1$ distinct states. The resulting system consists of $n+1$ coupled kinetic equations, which can be referred to as a ``multi‑species kinetic model" in which each component describes the evolution of a species associated with a fixed methylation level.

We divide $[M_l,M_r]$ into $n+1$ intervals, using grid points: 
$$
m_{-1/2}<M_l=m_0<m_{1/2}<m_{3/2}<\cdots<m_{n-1/2}<M_r=m_n<m_{n+1/2}.
$$
We discretize the advection term $\p_m\big (f(m-m_c) p(t,x,v,m) \big)$ using the standard finite volume method.
For $i=0,\cdots,n$, one can approximate $p(t,x,v,m)$ inside the interval $[m_{i-1/2},m_{i+1/2}]$ by the cell averages:
\begin{align*}
p_i(t,x,v)=\frac{1}{\Delta m_i}\int_{m_{i-1/2}}^{m_{i+1/2}}p(t,x,v,m)\,\d m,
\end{align*}
where $\Delta m_i=m_{i+1/2}-m_{i-1/2}$. Integrating both side of equation \eqref{PBKM1d} with respect to $m$ from $m_{i-1/2}$ to $m_{i+1/2}$, one gets
\begin{align}\label{discrete system}
    \p_t p_i + v\p_x p_i + \tfrac{F_{{i+1/2}}-F_{i-1/2}}{\Delta m_i} = \frac{1}{\Delta m_i}\int_{m_{i-1/2}}^{m_{i+1/2}} \Lambda \big(m-m_c\big)\big(\tfrac{1}{2}\int_{-1}^1 p \,\d v' - p\big)\,\d m,
\end{align}
where $F_{i\pm 1/2}$ are the fluxes at the cell interfaces $m_{i\pm 1/2}$.

One can approximate 
$$
\frac{1}{\Delta m_i}\int_{m_{i-1/2}}^{m_{i+1/2}} \Lambda (m-m_c(x))(\frac{1}{2}\int_{-1}^1 p \,\d v' - p)\,\d m \approx \Lambda_i (\frac{1}{2}\int_{-1}^1 p_i \,\d v' - p_i),
$$
where $\Lambda_i(x)=\Lambda (\frac{1}{2}(m_{i+1/2}+m_{i-1/2})-m_c)$. The upwind method can be employed to approximate the flux at the cell interfaces. Due to the property of $f(y)$ in \eqref{eq:f}, the signs of the sequence 
$$f(m_{-1/2}-m_c(x)),\, f(m_{1/2}-m_c(x)),\, \cdots,\, f(m_{n+1/2}-m_c(x))$$
change from positive to negative. Then we use the following approximations: when $f(m_{i+1/2}-m_c(x))\geq 0$, let
\begin{subequations}\label{upwindF}
\begin{align}\label{upwind f>0}
F_{i+1/2}= f(m_{i+1/2}-m_c(x))p_i(t,x,v),
\end{align}
and when $f(m_{i+1/2}-m_c(x))< 0$, 
let
\begin{align*}\label{upwind f<0}
F_{i+1/2}=f(m_{i+1/2}-m_c(x)) p_{i+1}(t,x,v).
\end{align*} 
\end{subequations}
 For all $x\in\R$, there exists an index $i_0\in\{0,1,\cdots,n\}$ such that the equilibrium methylation level $m_c(x)\in [m_{i_0-1/2},m_{i_0+1/2})$. Then $f(m_{i_0-1/2}-m_c(x))>0$ and $f(m_{i_0+1/2}-m_c(x))<0$. The upwind discretizations of convective terms of $m$ for $p_{i_0-1},\,p_{i_0+1}$ are respectively given by 
\begin{equation}\label{upwind discretization0}
\begin{aligned}
&    \frac{f(m_{i_0-1/2}-m_c)p_{i_0-1}(t,x,v)-f(m_{i_0-3/2}-m_c)p_{i_0-2}(t,x,v)}{\Delta m_{i_0-1}},\\  &\frac{f(m_{i_0+3/2}-m_c)p_{i_0+2}(t,x,v)-f(m_{i_0+1/2}-m_c)p_{i_0+1}(t,x,v)}{\Delta m_{i_0+1}},
\end{aligned}
\end{equation}
while the convective term in $m$ for $p_{i_0}$ is written as
\begin{align}\label{upwind discretization i}
  \frac{f(m_{i_0+1/2}-m_c)p_{i_0+1}(t,x,v) - f(m_{i_0-1/2}-m_c)p_{i_0-1}(t,x,v)}{\Delta m_{i_0}}.
\end{align}

In order to close the system,
the values of $p_{-1}(t,x,v)$ and $p_{n+1}(t,x,v)$ are set to be zero, i.e.,
\begin{align*}
    p_{-1}(t,x,v)=p_{n+1}(t,x,v)=0,
\end{align*}
according to the boundary conditions of the model \eqref{sys:kinetic-transport model}.

\paragraph{Conservation of mass at the discrete level} The total mass at time $t$ is defined as
\begin{align*}
    M(t) = \sum_{i=0}^{n} \int_{\mathbb{R}} \int_{-1}^1 p_i(t,x,v) \, \d v\, \d x \, \Delta m_i.
\end{align*}
Integrating equation \eqref{discrete system} over $(x,v) \in \mathbb{R} \times [-1,1]$, summing over $i = 0, \dots, n$, and multiplying each equation by $\Delta m_i$ yields
\begin{align}\label{mass term}
    \sum_{i=0}^{n}\int_{\mathbb{R}} \int_{-1}^1 \p_t p_i \, \d v\, \d x \, \Delta m_i = \frac{\d}{\d t} \sum_{i=0}^{n} \int_{\mathbb{R}} \int_{-1}^1 p_i \, \d v\, \d x \, \Delta m_i = \frac{\d}{\d t} M(t).
\end{align}
The spatial transport term is
\begin{align*}
\sum_{i=0}^{n} \int_{\mathbb{R}} \int_{-1}^1 v \p_x p_i \, \d v\, \d x \, \Delta m_i = \sum_{i=0}^{n} \int_{-1}^1 \left( \int_{\mathbb{R}} v \p_x p_i \, \d x \right) \d v \, \Delta m_i =0.
\end{align*}
Consider the advection term in \eqref{upwindF}-\eqref{upwind discretization i}, summing over $i = 0, \dots, n$  and using the boundary conditions, one gets
\begin{align}\label{advection term}
\sum_{i=0}^{n} \left(F_{i+\frac{1}{2}} - F_{i-\frac{1}{2}}\right) = F_{n+\frac{1}{2}} - F_{-\frac{1}{2}}=0.
\end{align}
Here the last equality follows from $p_{-1}=p_{n+1}=0$, which gives $F_{-\frac{1}{2}}=F_{n+\frac{1}{2}}=0$.

Integrating the collision term over $\mathbb{R} \times [-1,1]$, one has
\begin{align}\label{collision term}
\sum_{i=0}^{n} \int_{\mathbb{R}} \int_{-1}^1 \Lambda_i \left( \frac{1}{2}\int_{-1}^1 p_i(t,x,v') \d v' - p_i(t,x,v) \right) \d v\, \d x = 0.
\end{align}
Combining \eqref{mass term}-\eqref{collision term}, we obtain
$\frac{\d}{\d t} M(t) = 0$.
This proves that the total mass is conserved in time for the discrete system.

In the discrete model, each $p_i$ represents one species, while different species are coupled together by the discretized $m$-advection term. The signs in front of the discretized $m$-advection term indicate the transitions between different species. For example, in \eqref{upwind discretization0}, since $f(m_{i_0-1/2}-m_c)>0$ and $f(m_{i_0-3/2}-m_c)>0$, the term $f(m_{i_0-3/2}-m_c)p_{i_0-2}$ indicates that species $p_{i_0-2}$ transitions into $p_{i_0-1}$; since $f(m_{i_0+3/2}-m_c)<0$ and $f(m_{i_0+1/2}-m_c)<0$, the term $f(m_{i_0+3/2}-m_c)p_{i_0+2}$ means that species $p_{i_0+2}$ transitions into $p_{i_0+1}$. In \eqref{upwind discretization i}, since $f(m_{i_0-1/2}-m_c)>0$ and $f(m_{i_0+1/2}-m_c)<0$, the term $f(m_{i_0-1/2}-m_c)p_{i_0-1}$ (respectively, $f(m_{i_0+1/2}-m_c)p_{i_0+1}$) indicates that species $p_{i_0-1}$ (respectively, $p_{i_0+1}$) transitions into $p_{i_0}$.

\subsection{The interface problem for multi-species model and the scalings}
The physical domain is composed of two domains  $\Omega_l=(-\infty,0)$ and $\Omega_r=(0,+\infty)$. For simplicity of the notations, only uniform meshes are considered, i.e. $\Delta m_0=\cdots=\Delta m_{n}=\Delta m$, but the discretization and derivation can be easily generalized to nonuniform grids. For $i = 0,1,\dots,n$, let 
\begin{equation}\label{variable:Ril}
R_{i}^{l}(x)= \left|\frac{f(m_{i-1/2} - M_l)}{\Delta m}\right| (x\in\Omega_l), \,\,\,\,
R_{i}^{r}(x)= \left|\frac{f(m_{i+1/2} - M_r)}{\Delta m}\right| (x\in\Omega_r).
\end{equation}\noindent
By the definition of $R_{i}^{l}(x) (R_{i}^{r}(x))$ and for the convenience of the subsequent proof, we assume $R_{i}^{l}(x)(R_{i}^{r}(x))\geq c_0>0$. Additionally, to distinguish the definition of the tumbling term $\Lambda_i$ ($i = 0,1,\dots,n$) in the two domains $\Omega_l$ $(\Omega_r)$, we denote them as $\Lambda_i^l$ $(\Lambda_i^r)$ respectively.
 
By using $p_{-1}=p_{n+1}=0$, the discrete system inside $\Omega_l$ and $\Omega_r$ can be written as follows:
\begin{subequations}\label{sys:n-species}
\begin{align}\label{sys:n-species-l}
 x\in\Omega_l
  \begin{cases}
    \p_t p_0  + v\p_x p_0 - R_1^l p_1 = \Lambda_0^l \left(\tfrac{1}{2} \int_{-1}^1 p_0 \,\d v' - p_0\right),\\[2mm]
    \p_t p_i + v\p_x p_i + R_i^l p_i  - R_{i+1}^l p_{i+1} = \Lambda_i^l \left(\tfrac{1}{2} \int_{-1}^1 p_i \,\d v' - p_i\right),\\[2mm]
    \p_t p_{n} + v\p_x p_{n} + R_n^l p_{n} =  \Lambda_{n}^l \left(\tfrac{1}{2} \int_{-1}^1 p_{n} \,\d v' - p_{n}\right),
  \end{cases}
  \end{align}
\begin{align}\label{sys:n-species-r}
   x\in\Omega_r
 \begin{cases}
    \p_t p_0  + v\p_x p_0 + R_0^r p_0 =  \Lambda_0^r \left(\tfrac{1}{2}\int_{-1}^1 p_0 \,\d v' - p_0\right),\\[2mm]
    \p_t p_i + v\p_x p_i - R_{i-1}^r p_{i-1}  + R_i^r p_i = \Lambda_i^r \left(\tfrac{1}{2}\int_{-1}^1 p_i \,\d v' - p_i\right),\\[2mm]
    \p_t p_{n} + v\p_x p_{n} - R_{n-1}^r p_{n-1} =  \Lambda_{n}^r \left(\tfrac{1}{2}\int_{-1}^1 p_{n} \,\d v' - p_{n}\right),
  \end{cases}
\end{align}
\end{subequations}
where $i=1,2,\dots,n-1$. The densities for all species $p_i$ ($i=0,1,\dots,n$) are continuous at the interface $x = 0$. 

To identify a physically meaningful scaling for the interface problem, we start from two main assumptions as in \cite{Perthame-Tang-Vauchelet}. First, the intracellular adaptation dynamics drives the internal state $m$ toward its local equilibrium value $m_c(x)$ on a faster time scale than the spatial transport of the cells. We therefore introduce a small parameter $\varepsilon$ and assume that the adaptation dynamics occurs on a short time scale of $O(\varepsilon)$.
Second, the tumbling frequency $\Lambda$ is highly sensitive to the mismatch $m-m_c$. For example, in \textit{E. coli} chemotaxis, the tumbling response is often described by a Hill-type function with a large exponent \cite{Si2012}, so even a small deviation from equilibrium may lead to a large change in the tumbling frequency. When $m=m_c$, the tumbling rate is of order $O(1)$. When $m>m_c$, the tumbling frequency may become very large, of order $O(1/\varepsilon)$, whereas when $m<m_c$, the tumbling frequency may become very small.

This asymmetry becomes especially important near a sharp interface between two environments with very different external signals. Since the internal state $m$ cannot adapt instantaneously when a cell crosses the interface, the cell experiences a sudden mismatch between its current internal state and the new local equilibrium $m_c(x)$. Since $M_l<M_r$, cells moving from $\Omega_r$ to $\Omega_l$ experience a large positive mismatch and therefore tumble very frequently. In contrast, bacteria moving from $\Omega_l$ to $\Omega_r$ experience a mismatch that suppresses tumbling, allowing them to continue their run. This directional asymmetry is the microscopic origin of the biased scattering behavior observed at the interface.

We next nondimensionalize system \eqref{sys:n-species} according to the observations made above. For $s=0,1,\dots,n$, write
\begin{align*} 
  x =\hat{x}\, l_{\infty},
\quad 
  t= \hat{t}\, t_{\infty}, 
\quad 
  v = \hat{v}l_{\infty}/t_\infty,
\quad 
  R^{l,r}_s =\hat{R}^{l,r}_s/{t_{a}},
\quad 
  \Lambda_s^{l,r}= \hat\Lambda_s^{l,r}/t^{l,r}_{tum,s}. 
\end{align*}
Here, variables with hats denote nondimensional quantities, while quantities with the subscript $\infty$ denote characteristic reference scales with physical units. The parameters $t_a$ and $t_{tum,s}$ denote the characteristic adaptation time and the characteristic tumbling time for different species, respectively. We assume fast adaptation and a stiff tumbling response. Thus, $t_a$ is small, while the tumbling time $t^{l,r}_{tum,s}$ may have different asymptotic orders for different $s$, $l$, or $r$. This reflects the fact that different species may respond to the same external environment with very different tumbling frequencies.

When we say that a species is fully adapted or poorly adapted, we mean that its internal state is, respectively, closest to or farthest from the local equilibrium value $m_c$. In $\Omega_l$, cells adapt toward $M_l$. Therefore, the species $p_0$ is the best-adapted species and has an $O(1)$ tumbling time, thus we choose $t_{tum,0}^l=1$.  By contrast, the species $p_n$ is the most poorly adapted species in $\Omega_l$. Since it corresponds to a large positive mismatch, its tumbling frequency is very large and its characteristic tumbling time is small, $t_{tum,n}^l=\varepsilon$. The intermediate species $p_i$, with $1\le i\le n-1$, have intermediate tumbling time scales, but for simplicity we choose  $t_{tum,i}^l=\varepsilon$; then $\hat\Lambda_s^l$ can range from $O(1)$ to $O(\varepsilon)$ when $s$ changes from $n$ to $0$.

In $\Omega_r$, the local equilibrium is changed to $M_r$. Hence $p_n$ becomes the best-adapted species, and the species $p_0$ is now the most poorly adapted species. In this case, the mismatch suppresses tumbling, so $p_0$ has a very small tumbling frequency, or equivalently a very large characteristic tumbling time with $t_{tum,0}^r=\varepsilon^{-\beta_0}$ ($\beta_0>0$). The intermediate species $p_s$ ($s=1,\dots,n-1$) again have intermediate scalings with ${t_{tum,s}^r}=\varepsilon^{-\beta_s}$, where $\beta_s>0$. Moreover, $t_{tum,n}^r=1$ in this setting, since the fully adapted species has $O(1)$ tumbling frequency, hence $\hat\Lambda_n^r$ is $O(1)$.

Incorporating the above-mentioned fast adaptation and asymmetric tumbling scalings into the system \eqref{sys:n-species}, we obtain the following nondimensionalized system (by dropping the hats):
\begin{subequations}\label{sys:n-species-scaling}
\begin{align}\label{sys:n-species-0}
 x\in\Omega_l
  \begin{cases}
    \p_t p^{\eps}_0  + v\p_x p^{\eps}_0 - \frac{1}{\eps}R_1^l p^{\eps}_1 = \Lambda_0^l \left(\tfrac{1}{2} \int_{-1}^1 p^{\eps}_0 \,\d v' - p^{\eps}_0\right),\\[2mm]
    \p_t p^{\eps}_i + v\p_x p^{\eps}_i + \frac{1}{\eps}R_i^l p^{\eps}_i  - \frac{1}{\eps}R_{i+1}^l p^{\eps}_{i+1} = \frac{1}{\eps}\Lambda_i^l \left(\tfrac{1}{2} \int_{-1}^1 p^{\eps}_i \,\d v' - p^{\eps}_i\right),\\[2mm]
    \p_t p^{\eps}_{n} + v\p_x p^{\eps}_{n} + \frac{1}{\eps}R_n^l p^{\eps}_{n} = \tfrac{1}{\eps} \Lambda_{n}^l \left(\tfrac{1}{2} \int_{-1}^1 p^{\eps}_{n} \,\d v' - p^{\eps}_{n}\right),
  \end{cases}
\end{align}
\begin{align}\label{sys:n-species-N}
   x\in\Omega_r
 \begin{cases}
    \p_t p^{\eps}_0  + v\p_x p^{\eps}_0 + \frac{1}{\eps}R_0^r p^{\eps}_0 =  \eps^{\beta_0} \Lambda_0^r \left(\tfrac{1}{2}\int_{-1}^1 p^{\eps}_0 \,\d v' - p^{\eps}_0\right),\\[2mm]
    \p_t p^{\eps}_i + v\p_x p^{\eps}_i - \frac{1}{\eps}R_{i-1}^r p^{\eps}_{i-1}  + \frac{1}{\eps}R_i^r p^{\eps}_i = \eps^{\beta_i}\Lambda_i^r \left(\tfrac{1}{2}\int_{-1}^1 p^{\eps}_i \,\d v' - p^{\eps}_i\right),\\[2mm]
    \p_t p^{\eps}_{n} + v\p_x p^{\eps}_{n} - \frac{1}{\eps}R_{n-1}^r p^{\eps}_{n-1} =  \Lambda_{n}^r \left(\tfrac{1}{2}\int_{-1}^1 p^{\eps}_{n} \,\d v' - p^{\eps}_{n}\right),
  \end{cases}
\end{align}
\end{subequations}
where $\beta_0, \cdots, \beta_{n-1}>0$ and  $\Lambda_1^{l}, \cdots, \Lambda_n^{l}$, $\Lambda_0^{r}, \cdots, \Lambda_{n-1}^r$ can range from $O(1)$ to $O(\varepsilon)$.

\paragraph{The interface conditions} The densities for all species $p_i^{\eps}$ ($i=0,1,\dots,n$) are continuous at the interface $x = 0$ such that
 \begin{align}\label{BC:continuous interface-n-conti}
p_i^{\eps} (t,0^-,v) = p_i^{\eps} (t,0^+,v), \quad i=0,1,\dots,n.
\end{align}

\section{Asymptotic Analysis}\label{sec:Main results}
In this section, we formally derive the interior equations in $\Omega_l$, $\Omega_r$ together with their interface conditions via an asymptotic analysis of the multi-species model \eqref{sys:n-species-scaling} in the limit $\eps \to 0$. We focus on the interface layer and assume that the initial data are well prepared so that there will be no initial layer generated. This assumption can be removed by analyzing initial and initial-boundary layers to make the following formal analysis rigorous. We nevertheless opt for a formal analysis, since our focus is on the derivation of the interface condition rather than the technical justification of the expansion. 

\subsection{Statement of the limiting system}
    Let $p_i^{\eps}(t,x,v)$\,\,$(i=0,1,\dots,n)$ be the solutions to the system \eqref{sys:n-species-scaling}-\eqref{BC:continuous interface-n-conti}. As $\eps\rightarrow 0$, $p^{\eps}(t,x,v) = \sum_{i=0}^n p_i^{\eps}$ converges to $P(t,x,v)$ that satisfies the following kinetic model in each subdomain:
   \begin{equation}
   \begin{aligned} \label{sys:conti-limit-n}
     \p_t P  + v\p_x P &= \Lambda_0^l (\tfrac{1}{2} \int_{-1}^1 P \, \d v' - P), \quad x\in \Omega_l,\\
      \p_t P  + v\p_x P &= \Lambda_n^r (\tfrac{1}{2} \int_{-1}^1 P \, \d v' - P), \quad x\in \Omega_r.
   \end{aligned}
   \end{equation}
 Furthermore, there exists an operator $\mathscr{A}$, depending on $\Lambda_i^l$  and  $R_i^l$, $i = 1, \cdots, n$,  such that the interface condition of the limiting kinetic model satisfies
\begin{align} \label{cond:transmission}
     v P(t,0^-, -v) = \mathscr{A} [\tau_-P (t,0^+, \cdot)],
\qquad  v \in (-1, 1),
   \end{align}
where $\tau_-$ is the reflection operator $(\tau_- f)(v) = f(-v)$ such that $\mathscr{A}$, as a nonlocal operator, acts on $P(t, 0^+, -v')$ with $v' \in (-1, 1)$. 
\begin{remark}
Note that even when the respective tumbling rates $\Lambda^l_0$ and $\Lambda^r_n$ in $\Omega_l$ and $\Omega_r$ are the same, the density function $P$ still experiences a jump at the interface given by~\eqref{cond:transmission}. This is induced by the sharp transition of $M_l$ to $M_r$ across the interface with infinitesimal width.
\end{remark}

\subsection{Formal derivation} The derivation consists of two parts: the interior equations and the interface condition. By the standard comparison principle, $p_i^\eps(t,x,v)$ is nonnegative for all $i, \eps$, if their initial data are nonnegative. Moreover, the continuity of $p^\eps_i$ at the interface $x=0$ guarantees that the mass conservation holds:
\begin{align*}
  \sum_{i=0}^n\int_\R \int_{-1}^1 p_i^\eps (t, x, v) \,\d x \,\d v = \sum_{i=0}^n\int_\R \int_{-1}^1 p_i^\eps (0, x, v) \,\d x \,\d v < \infty.  
\end{align*}
Therefore, there exists a probability measure $P \geq 0$ such that 
\begin{align*}
  \text{$\sum_{i=0}^n p_i^\eps \to P$ as $\eps \to 0$ weakly as measures up to a subsequence.}
\end{align*}

\noindent \textit{Step 1. Derivation of the interior equations.} We start with deriving the limiting equation in the left-side domain $\Omega_l = (-\infty, 0)$. Specifically, we show inductively that $p_i^\eps \to 0$ in weighted norms for any $i = 1, 2, \cdots, n$, due to the dissipation of both the adaptation and the tumbling terms. 
We will prove that $p_i^\eps \to 0$ in $L^1(w \d x \d v)$ through the energy estimate, where the weight function is written as
\begin{align*}
  w(x)=x^2 e^x, \qquad x \in \Omega_l = (-\infty, 0).
\end{align*}
The particular choice of $w(x)$ is to remove the interface and boundary conditions at $x = 0, -\infty$. Since we only have the $L^1$-boundedness of $p_i^\eps$ (instead of pointwise boundedness), we further introduce two cutoff functions $\phi, \psi \in C^\infty(-\infty, 0)$ such that 
\begin{align*}
  \phi(x)
= \begin{cases}
   1, & x \in (-\infty, -2], \\[3pt]
   0, & x \in (-1, 0],
  \end{cases}
\qquad
  \psi(x)
= \begin{cases}
   1, & x \in [-1, 0],  \\[3pt]
   0, & x \in (-\infty, -2]. 
  \end{cases}
\end{align*}
For any $\delta > 0$, let 
\begin{align*}
  \phi_\delta = \phi(x/\delta),
\qquad
  \psi_\delta = \psi (\delta x), 
\qquad
  \text{for $x \in \Omega_l$}. 
\end{align*}
More explicitly, 
\begin{align*}
  \phi_\delta(x)
= \begin{cases}
   0, & x \in [-\delta, 0], \\[3pt]
   1, & x \in (-\infty, -2\delta],
  \end{cases}
\qquad
  \psi_\delta(x)
= \begin{cases}
   1, & x \in (-1/\delta, 0], \\[3pt]
   0, & x \in (-\infty, -2/\delta].
  \end{cases}
\end{align*}
To derive the full energy estimates for the full system, we start with the energy estimate for $p_n^\eps$. Multiply the $p_n^\eps$-equation by $w \phi_\delta \psi_\delta$ and integrate over $\Omega_l\times(-1,1)$. The transport term gives
\begin{align*}
\int_{\Omega_l}\int_{-1}^1 w v \phi_\delta \psi_\delta\p_x p_n^\eps\,\d v\,\d x
=
- \int_{\Omega_l}\int_{-1}^1 v (w \phi_\delta \psi_\delta)' p_n^\eps\,\d v\,\d x,
\end{align*}
where the boundary terms vanish due to the supports of the cutoff functions $\psi_\delta$ and $\phi_\delta$. For each $x \in (-\infty, 0)$, by the definitions of $\phi_\delta$ and $\psi_\delta$, we have
\begin{align*}
  (w \phi_\delta \psi_\delta)'(x)
= w'(x) (\phi_\delta \psi_\delta)(x)
  + w(x) (\phi_\delta \psi_\delta)'(x)
\to 
  w'(x). 
\end{align*}
Moreover, 
\begin{align*}
  |(w \phi_\delta \psi_\delta)'(x)|
&\leq
  |w'| \phi_\delta \psi_\delta
  + w |\phi'_\delta| \psi_\delta
  + w \phi_\delta |\psi'_\delta|
\\
& 
\leq
  |w'| \phi_\delta \psi_\delta
  + x e^x \frac{x}{\delta} |\phi'(x/\delta)| \psi_\delta
  + \delta w \phi_\delta |\psi'(\delta x)|
\\
& \leq 
 C( |w'(x)| + |x| e^x + w(x))  \in L^1(-\infty, 0),
\end{align*}
since $\frac{x}{\delta} \phi'(x/\delta) \in L^\infty(-\infty, 0)$ due to the fact that $\phi'$ is compactly supported on $[-2, -1]$. 
By the Lebesgue dominated convergence theorem, 
\begin{align*}
  \int_{\Omega_l}\int_{-1}^1 w v \phi_\delta \psi_\delta\p_x p_n^\eps\,\d v\,\d x
\to 
  - \int_{\Omega_l}\int_{-1}^1 w' v p_n^\eps\,\d v\,\d x
\qquad
  \text{as $\delta \to 0$.}
\end{align*}
Hence, the weighted $L^1$-energy estimate satisfies
\begin{align} \label{est:L1-pn}
\frac{\d}{\d t} \|w p_n^\eps\|_{L^1_{x,v}}
+\frac{R_n^l}{\eps}\|w p_n^\eps\|_{L^1_{x,v}}
= \int_{\Omega_l}\int_{-1}^1 v w'p_n^\eps\,\d v\,\d x,
\end{align}
where we have used that $\int_{-1}^1 \left(\tfrac{1}{2}\int_{-1}^1 p^{\eps}_{n} \,\d v' - p^{\eps}_{n}\right) \,\d v = 0$. 

Since $\|p_n^\eps\|_{L^1_{x,v}}$ is bounded uniformly in $\eps$, for any $c_1 > 0$ and $\eps$ small enough such that $c_1 \eps < 1$, 
\begin{align*}
\int_{\Omega_l}\int_{-1}^1 |v|\,|w'|\,|p_n^\eps| \,\d v\,\d x
&=
\int_{-\infty}^{-c_1\eps} \int_{-1}^1 |v||w'|\,|p_n^\eps|\,\d v\,\d x
+  \int_{-c_1 \eps}^0 \int_{-1}^1 |v||w'|\,|p_n^\eps|\,\d v\,\d x
\\
& \hspace{-1.8cm}
\leq
 \int_{-\infty}^{-c_1 \eps} \int_{-1}^1 |v|(x^2 + 2|x|)e^x\,|p_n^\eps|\,\d v\,\d x
+ 3 \int_{-c_1 \eps}^0 \int_{-1}^1 |v||x| e^x\,|p_n^\eps|\,\d v\,\d x
\\
& \hspace{-1.8cm}
\leq
\vpran{1+\frac{2}{c_1 \eps}} \|w p_n^\eps\|_{L^1_{x,v}}
+ C \eps,
\end{align*}
where $C > 0$ is independent of $\varepsilon$. If $c_1 > \frac{8}{R^l_n}$ and $\frac{\eps}{R^l_n} < \frac{1}{4}$ such that $1+\frac{2}{c_1 \eps} < \frac{R^l_n}{2\varepsilon}$, then \begin{align} \label{bound:w-deriv-p-n}
\int_{\Omega_l}\int_{-1}^1 |v|\,|w'|\,|p_n^\eps|\,\d v\,\d x
\le
\frac{R_n^l}{2\eps}\|w p_n^\eps\|_{L^1_{x,v}}
+ C \eps. 
\end{align}
Then \eqref{est:L1-pn} becomes
\begin{align*}
\frac{\d}{\d t}\|w p_n^\eps\|_{L^1_{x,v}}
+\frac{R_n^l}{2\eps}\|w p_n^\eps\|_{L^1_{x,v}}
\le
 C \eps,
\end{align*}
which gives
\begin{align}\label{eq:weighted norm pn}
\|w p_n^\eps(t)\|_{L^1_{x,v}}
\le
e^{- R_n^l \frac{t}{2\eps}}\|w p_n^\eps(0)\|_{L^1_{x,v}}
+ C \eps^2. 
\end{align} 

We apply similar estimates to bound the weighted $L^1$-norms of $p_{i}^\eps$ with $i=1,2,\dots, n-1$. 
In particular, by multiplying by \(w\), using the cutoff functions, integrating in $x,v$, and passing to the limit, the boundary term from the transport part vanishes, and we get
\begin{align*}
\frac{\d}{\d t}\|w p_i^\eps\|_{L^1_{x,v}}
+\frac{R_i^l}{\eps}\|w p_i^\eps\|_{L^1_{x,v}}
\le
\frac{R_{i+1}^l}{\eps}
\|w p_{i+1}^\eps\|_{L^1_{x,v}}
+
\int_{\Omega_l}\int_{-1}^1 |v|\,|w'|\,|p_i^\eps|\,\d v\,\d x.
\end{align*}
By a similar estimate as in \eqref{bound:w-deriv-p-n}, for $\eps$ small enough,
\begin{align*}
\int_{\Omega_l}\int_{-1}^1 |v|\,|w'|\,|p_i^\eps|\,\d v\,\d x
\le
\frac{R_i^l}{2\eps}\|w p_i^\eps\|_{L^1_{x,v}}
+
C \eps.
\end{align*}

Combining the above estimates, we have
\begin{align*}
\frac{\d}{\d t}\|w p_i^\eps\|_{L^1_{x,v}}
+\frac{R_i^l}{2\eps}\|w p_i^\eps\|_{L^1_{x,v}}
\le
C \frac{e^{-c_0 \frac{t}{\eps}}}{\eps} 
+ C \eps,
\qquad
c_0 < \tfrac 12 \min\{R^l_i\}_{i=1}^n,
\end{align*}
which, by Gronwall, leads to
\begin{align*}
\|w p_i^\eps(t)\|_{L^1_{x,v}}
\le
C_1 e^{-c_0 t/\eps}
+
C_2\eps^2, 
\end{align*}
where $C_1, C_2$ are independent of $t, \eps$ but they could depend on $n$.
Therefore, for the number of species $n$ fixed, outside the initial layer of $t = O(\varepsilon)$, each $p^\eps_i$ satisfies that $p_{i}^\eps \to 0$ for each $i \geq 1$. We are only left with the limiting equation for $p_0^\eps$. Denote the limit of $p_0^\eps$ as $p^l_0$. By passing $\eps \to 0$ and using the previous estimate of $p^\eps_1$ for $t \gg O(\eps)$, we obtain that 
\begin{align*}
\p_t p_0^l + v \p_x p_0^l
= \Lambda_0^l \Bigl( \tfrac{1}{2}\int_{-1}^1 p_0^l \,\d v' - p_0^l \Bigr),\quad x\in \Omega_l.
\end{align*}

Similarly, by inductive estimates starting from $i = 0$, we obtain the limiting equation for $p_n^\eps$ in the domain $\Omega_r$ as
 \begin{align} \label{eq:interior-pn-r}
    \p_t p_n^r + v \p_x p_n^r  &= \Lambda_n^r  (\tfrac{1}{2}\int_{-1}^1 p_n^r \d v' - p_n^r),\quad x\in \Omega_r.
\end{align}
Unify the notation by denoting
\begin{align*}
  P(t,x,v)
= \begin{cases}
   p^l_0 (t,x,v), & \text{in $\Omega_l$,} \\
   p^r_n (t,x,v), & \text{in $\Omega_r$}.
  \end{cases}
\end{align*}
Then $P$ satisfies the model~\eqref{sys:conti-limit-n} in the interior of $\Omega_l$ and $\Omega_r$.

\medskip

\noindent \textit{Step 2. Derivation of the interface condition.}
The interface condition between the domains $\Omega_l$ and $\Omega_r$ is derived by an interface layer analysis at $x=0$. To this end, we derive two boundary layer equations by rescaling $x$ near $x=0$. In particular, in the left domain $\Omega_l$, we let
\begin{align*}
   y= - \tfrac{x}{\eps} \,\,\text{for}\,\, x\in\Omega_l
\qquad \text{and} \qquad
  \tilde p_{i}^\eps (t,y, v)=p_i^{\eps}(t, -\eps y, -v), 
\end{align*}
where we have followed the convention of taking $v>0$ as the incoming data at the boundary $y=0$. Denote the leading order of each $\tilde p_{k}^\eps$ as $\tilde p_{k}$ for $k = 0, 1, \cdots, n$. Then they satisfy the half-space boundary layer equation on the left side:
\begin{equation} \label{eq:BL-L}
 \Omega_l: 
 \begin{cases}
    v\p_y \tilde p_0 -  R_1^l \tilde p_1 = 0, 
\\[1mm]
    v\p_y \tilde p_i +  R_i^l \tilde p_i  -  R_{i+1}^l \tilde p_{i+1} = \Lambda_{i}^l (\tfrac{1}{2} \int_{-1}^1 \tilde p_{i} \,\d v' - \tilde p_{i}), 
\\[1mm]
   v\p_y \tilde p_{n} +  R_n^l \tilde p_{n} =  \Lambda_{n}^l (\tfrac{1}{2} \int_{-1}^1 \tilde p_{n} \,\d v' - \tilde p_{n}),  
  \end{cases}
\end{equation}
with $(y, v) \in (0, \infty) \times (-1, 1)$. The incoming data is to be prescribed at $y=0$. To be consistent with $\Omega_l$, in the right domain $\Omega_r$, we let
\begin{align*}
  z = \tfrac{x}{\eps} \quad \text{for $x \in \Omega_r$}
\qquad \text{and} \qquad
  \hat p^\eps_i (t, z, v)
= p^\eps_i (t, \eps z, -v) \quad \text{in $\Omega_r$.}
\end{align*}
Then the leading orders of $\widehat p^\eps_0, \cdots, \widehat p^\eps_n$ satisfy
\begin{equation} \label{eq:BL-R}
 \Omega_r: 
 \begin{cases}
  v\p_z \widehat p_0 -  R_0^r \widehat p_0 = 0,
  \\[1mm]
  v\p_z \widehat p_i +  R_{i-1}^r \widehat p_{i-1}  -  R_i^r \widehat p_i = 0, \\[1mm]
  v\p_z \widehat p_{n} +  R_{n-1}^r \widehat p_{n-1} =  0, 
  \end{cases}
\end{equation}
with $(z, v) \in (0, \infty) \times (-1,1)$.
The continuity~condition holds at the interface $y=0$ or~$z=0$:
\begin{align} \label{eq:conti-interf}
  \widetilde p_i (0, v) = \widehat p_i (0, v), 
\qquad
  v \in (-1, 1),
\quad
  i = 0, 1, \cdots, n,
\end{align}
where we have suppressed the time variable. 

\smallskip

\noindent \underline{\textbf{Analysis of~\eqref{eq:BL-R}}.} \,
Integrating the first equation in \eqref{eq:BL-R} gives, for
$v\neq 0$,
\[
\widehat p_0(z,v)=e^{R_0^r z/v}\widehat p_0(0,v).
\]
Boundedness requires $\widehat p_0(0,v)=0$ for $v>0$. Hence,
\begin{equation}\label{soln:hat-p-0-posit}
\widehat p_0(z,v)=
\begin{cases}
e^{R_0^r z/v}\widehat p_0(0,v), & v<0,\\
0, & v>0.
\end{cases}
\end{equation}

Substituting $\widehat p_0 (z,v)$ into the other equations in~\eqref{eq:BL-R} and solving inductively gives 
\begin{align*}
\widehat{p}_i(z, v) &= e^{\frac{R_i^r}{v} z} \widehat{p}_i(0, v) - \frac{R_{i-1}^r}{v} \int_0^z e^{\frac{R_i^r}{v} (z - s)} \widehat{p}_{i-1}(s, v) \, \d s, \quad i = 1, \ldots, n-1, \\
\widehat{p}_n(z, v) &= \widehat{p}_n(0, v) - \frac{R_{n-1}^r}{v} \int_0^z \widehat{p}_{n-1}(s, v) \, \d s.
\end{align*}
Two observations are in order. First, given the boundedness of each $\widehat p_i$ for $i = 0, \cdots, n-1$ and $\widehat p_0 = 0$ for $v > 0$, we can derive inductively that 
\begin{align} \label{eq:pi-v-posit}
   \widehat p_i (z, v) = 0, 
\qquad 
  \text{for all $v > 0$, $z > 0$, $i = 0, 1, \cdots, n-1$.} 
\end{align}
As a consequence, 
\begin{align} \label{eq:pnr-v-posit}
  \widehat p_n (z, v) = \widehat p_n (0, v), 
\qquad
  v > 0. 
\end{align}
Meanwhile, for $v < 0$, an inductive estimate using Lemma SM1.1 (see Appendix) gives 
\begin{align} \label{lim:p-i-r-neg}
   \widehat p_i (\infty, v) = 0, 
\qquad 
  \text{for all $v < 0$ and $i = 0, 1, \cdots, n-1$.} 
\end{align}
Therefore,
\begin{align} \label{lim:p-i-r}
   \widehat p_i (\infty, v) = 0, 
\qquad 
  \text{for all $v \in (-1, 1)$ and $i = 0, 1, \cdots, n-1$.} 
\end{align}
Second, by adding up all the equations in system~\eqref{eq:BL-R}, one derives the conservation law
\begin{align*}
  v \p_z \sum_{i=0}^n \widehat p_i = 0. 
\end{align*}
Therefore, by~\eqref{eq:pi-v-posit} and~\eqref{lim:p-i-r},
\begin{align} \label{eq:pn-inf}
  \widehat p_n (\infty, v)
= \sum_{i=0}^n \widehat p_i (\infty, v)
= \sum_{i=0}^n \widehat p_i (0, v)
= \begin{cases}
   \widehat p_n (0, v), & v > 0, \\[2pt]
   \sum_{i=0}^n \widehat p_i (0, v), & v < 0.
  \end{cases}
\end{align}
Note that in the interior of $\Omega_r$, only species $n$ survives and equation~\eqref{eq:pn-inf} provides the boundary condition for the interior equation~\eqref{eq:interior-pn-r}. 

\medskip 

\noindent \underline{\textbf{Analysis of~\eqref{eq:BL-L}}.} \,
For the left side domain $\Omega_l$,  by the interface continuity condition \eqref{eq:conti-interf}, the boundary conditions in $\Omega_l$ can be expressed in terms of the solution in $\Omega_r$. This establishes a half-space problem in $\Omega_l$:
\begin{align}\label{half-space-nspecies}
 \begin{cases}{}
   v \p_y \widetilde p_n + R_n^l\widetilde p_n = \Lambda_n^l (\tfrac{1}{2}\int_{-1}^1 \widetilde p_n \d v' - \widetilde p_n), 
   \\
   \widetilde p_n (0,v)= \widehat p_n(0, v),\quad v>0,
   \\
   \text{$\widetilde p_n$ is bounded.}
  \end{cases}
  \end{align}
Denote 
\begin{align*}
   \mathcal{L}_R f:= \Lambda_n^l \vpran{\tfrac{1}{2}\int_{-1}^1 f(v') \d v' - f} - R_n^l f. 
\end{align*}
Then the half-space equation has a bounded solution, which is uniquely determined by the incoming data $\widehat p_n(0, v)$ \cite{CGS-1988, LLS-2017a, LLS-2017b}. Therefore, there exists an operator $\mathscr{A}_n$ such that
\begin{align} \label{def:A-1}
  \widetilde p_n (y, v)
= \mathscr{A}_n [\widehat p_n (0, \cdot)]
= \mathscr{A}_n [\widehat p_n (\infty, \cdot)],
\qquad
(y, v) \in [0, \infty) \times (-1, 1), 
\end{align}
where $\mathscr{A}_n$ only acts on the incoming data where $v > 0$. Moreover,
\begin{align*}
  \widetilde p_n (\infty, v) \in \text{Null}(\mathcal{L}_R),
\end{align*}
with an exponential convergence rate~\cite{Cercignani-1988}. Since $\mathcal{L}_R$ is strictly dissipative, its null space is trivial. Therefore, 
\begin{align*}
  \widetilde p_n (y, v) \to 0
\quad \text{exponentially as $y \to \infty$.}
\end{align*}

Similarly, for each $i = 1, \cdots, n-1$, $\widetilde p_i$ satisfies a forced half-space equation:
\begin{align*}
 \begin{cases}{}
   v \p_y \widetilde p_i + R_i^l\widetilde p_i = \Lambda_i^l (\tfrac{1}{2}\int_{-1}^1 \widetilde p_i \d v' - \widetilde p_i) + R_{i+1}^l \widetilde p_{i+1}, 
   \\
   \widetilde p_i (0,v)= 0,\quad v>0,
   \\
   \text{$\widetilde p_i$ is bounded.}
  \end{cases}
\end{align*}
Note that the zero incoming data condition for $\widetilde p_i$ follows from the continuity at the boundary and~\eqref{eq:pi-v-posit}. Moreover, since $\widetilde p_n (y, v)$ decays exponentially as $y \to \infty$, we derive inductively that $\widetilde p_i (y, v)$ decays exponentially as $y \to \infty$ \cite{CGS-1988}. Furthermore, by induction again, since $\widetilde p_i$ depends on $\widetilde p_{i+1}$ and by~\eqref{def:A-1}, $\widetilde p_n$ depends on $\widehat p_n(\infty, v)$, we can denote each solution as 
\begin{align*}
\hspace{2cm}
  \widetilde p_i(y, v) = \mathscr{A}_i[\widehat p_n(\infty, \cdot)], 
\qquad
   y \in [0, \infty), \,\, v \in (-1, 1), \,\, i = 1, 2, \cdots, n-1. 
\end{align*}

Summing all the equations in~\eqref{eq:BL-L} and using the notations in~\eqref{def:A-1}, we obtain
\begin{align*}
v \p_y \sum\limits_{i=0}^{n}  \widetilde p_i (y,v) 
= \sum_{i=1}^n \Lambda_i^l \Big(\underbrace{\tfrac{1}{2}\int_{-1}^1 \widetilde p_i (y,v') \d v' - \widetilde p_i (y,v)}_{=:\mathcal{L} (\tilde p_i)}\Big)
= \sum_{i=1}^n \Lambda_i^l \mathcal{L}\vpran{\mathscr{A}_i [\widehat p_n (\infty, \cdot)]}.
 \end{align*}
Integrating $y$ from $0$ to infinity and using~\eqref{eq:pn-inf}, 
\begin{align} \label{def:A}
 v \widetilde p_0 (\infty,v) 
& = \sum\limits_{i=0}^{n} v \widetilde p_i (0,v) 
   + \sum_{i=1}^n \Lambda_i^l \int_0^{\infty}
   \mathcal{L}\vpran{\mathscr{A}_i [\widehat p_n (\infty, \cdot)]} (y, v) \d y \nn
\\
& = v \widehat p_n (\infty, v)
    + \sum_{i=1}^n\Lambda_i^l \int_0^{\infty}
   \mathcal{L}\vpran{\mathscr{A}_i [\widehat p_n (\infty, \cdot)]} (y, v) \d y 
 =: \mathscr{A} [\widehat  p_n (\infty,\cdot)].
\end{align}
Note that $y = \infty$ and $z = \infty$ correspond to the boundary data for the interior solutions, that is, 
\begin{align} \label{data:int-bdry}
 v \widetilde p_0(t, \infty, v) = v P (t, 0^-, -v), 
\qquad
  v \widehat p_n(t, \infty, v) = v P (t, 0^+, -v). 
\end{align}
Changing back to the physical velocity by reflection $v \to -v$, we obtain the desired transmission relation~\eqref{cond:transmission} at the interface.

\subsection{An example of a three-species system with two velocities}\label{subsec:three-species}

The operator $\mathscr{A}$ is defined implicitly through the solution
of a half-space problem. To make the interface law more explicit,
we give an example of a two-velocity model in which particles move only
in the forward and backward directions. This reduces the
integro-differential transport equation to two coupled PDE systems
in one spatial dimension. Such reductions are analogous to
the two-flux approximations commonly used in radiative
transfer~\cite{Modest2013}.

The interface operator is obtained by solving a coupled $2n+2$ linear ODE system and then eliminating the intermediate species. Although this procedure is constructive, the resulting transition matrix has a highly recursive and cumbersome explicit form. Therefore, instead of presenting the full $(n+1)$-species transition matrix, we use the three-species case as a representative example to illustrate the structure of the interface operator. This example already contains the essential features of the general multi-species setting.  

We consider the two-velocity model with $v\in\{-1,1\}$. Any common
speed magnitude $v_0>0$ can be reduced to this case by rescaling
the spatial coordinate. Let $p_i^{\pm,\eps}(t,x):=p_i^\eps(t,x,\pm1)$ for $i=0,1,2$.
The three-species, two-velocity system corresponding to
\eqref{sys:n-species-0} and \eqref{sys:n-species-N} in
$\Omega_l$ and $\Omega_r$, respectively, is given by
\begin{subequations}\label{sys:M-3-species}
\begin{align}\label{sys:M=m0-3-species}
  \Omega_l: \,\, \begin{cases}
    \p_t p^{\pm,\eps}_0 \pm \p_x p^{\pm,\eps}_0 - \frac{1}{\eps}R_1^l p^{\pm,\eps}_1 = \frac{1}{2}\Lambda_0^l (p^{\mp,\eps}_0 - p^{\pm,\eps}_0),\\[2mm]
    \p_t p^{\pm,\eps}_1 \pm \p_x p^{\pm,\eps}_1 + \frac{1}{\eps}R_1^l p^{\pm,\eps}_1 - \frac{1}{\eps} R_{2}^l p^{\pm,\eps}_{2} = \frac{1}{2\eps}\Lambda_1^l (p^{\mp,\eps}_1 - p^{\pm,\eps}_1),\\[2mm]
    \p_t p^{\pm,\eps}_{2} \pm \p_x p^{\pm,\eps}_{2} + \frac{1}{\eps}R_2^l p^{\pm,\eps}_{2} = \frac{1}{2\eps}\Lambda_{2}^l (p^{\mp,\eps}_{2} - p^{\pm,\eps}_{2}),
  \end{cases}
\end{align}
\begin{align}\label{sys:M=mn-3-species}
  \Omega_{r}: \,\,  \begin{cases}
    \p_t p^{\pm,\eps}_0 \pm \p_x p^{\pm,\eps}_0 + \frac{1}{\eps} R_0^{r} p^{\pm,\eps}_0 = \frac{1}{2}\eps^{\beta_0}\Lambda_0^r (p^{\mp,\eps}_0 - p^{\pm,\eps}_0),\\[2mm]
    \p_t p^{\pm,\eps}_1 \pm \p_x p^{\pm,\eps}_1 - \frac{1}{\eps}R_{0}^{r} p^{\pm,\eps}_{0} + \frac{1}{\eps}R_{1}^r p^{\pm,\eps}_{1} = \frac{1}{2}\eps^{\beta_1}\Lambda_1^r (p^{\mp,\eps}_1 - p^{\pm,\eps}_1),\\[2mm]
    \p_t p^{\pm,\eps}_2 \pm \p_x p^{\pm,\eps}_{2} - \frac{1}{\eps}R_{1}^r p^{\pm,\eps}_{1} =  \frac{1}{2}\Lambda_2^r (p^{\mp,\eps}_{2} - p^{\pm,\eps}_2),
  \end{cases}
\end{align}
\end{subequations}
where $\beta_0$, $\beta_1>0$.
Moreover, the densities for each species are continuous at the interface such that, for $i=0,1,2$,
\begin{align}\label{BC:connecting condition-3s}
  p_i^{\pm,\eps} (t,0^-) = p_i^{\pm,\eps} (t,0^+).
\end{align}
\noindent \textit{Statement of the limiting system.}
As $\eps\rightarrow 0$, the solution $p^{\pm,\eps} = \sum_{i=0}^2 p_i^{\pm,\eps}$ converges to $p^{\pm} = \sum_{i=0}^2 p_i^{\pm}$ satisfying
\begin{subequations}\label{eq:p-0-two-v-3s}
\begin{align}\label{eq:p-0-l-two-v-3s}
\Omega_l:
\begin{cases}
\p_t p^{+} + \p_x p^{+}  = \frac{1}{2}\Lambda_0^l (p^{-} - p^{+}),\\
\p_t p^{-} - \p_x p^{-}  = \frac{1}{2}\Lambda_0^l (p^{+} - p^{-}),
\end{cases}\\
\Omega_r:
\begin{cases}\label{eq:p-0-r-two-v-3s}
\p_t p^{+} + \p_x p^{+}  = \frac{1}{2}\Lambda_2^r (p^{-} - p^{+}),\\
\p_t p^{-} - \p_x p^{-}  = \frac{1}{2}\Lambda_2^r (p^{+} - p^{-}).
\end{cases}
\end{align}
\end{subequations}
Moreover, there exists a $2\times 2$ matrix $A_2$, depending on $\Lambda_1^l$, $\Lambda_2^l$, $R_1^l$ and $R_2^l$, such that the interface condition for $p^\pm$ at $0$ is given by
\begin{align*}
     \left(p^{+} (t,0^+), p^{-} (t,0^+)\right)^\top |_{\Omega_r}
= A_2 \,\, \left(p^{+} (t,0^-), p^{-} (t,0^-)\right)^\top|_{\Omega_l},
\end{align*}
where
\begin{align}\label{matrix:A-3 species}
A_2=\left(\begin{array}{cc}
1 &
\dfrac{(E+\lambda_1)D+R_2^l F}{(E+\lambda_1)(k-1)D-R_2^l F}\\[3mm]
0 &
\dfrac{k(E+\lambda_1)D}{(E+\lambda_1)(k-1)D-R_2^l F}
\end{array}\right),
\end{align}
with
\begin{align*}
&E=R_1^l+\tfrac{\Lambda_1^l}{2},\qquad
D=(R_1^l)^2+R_1^l\Lambda_1^l-\lambda_2^2,\qquad
k=\tfrac{\Lambda_2^l+2R_2^l+2\lambda_2}{\Lambda_2^l},\\
&F=
(E+\lambda_1)\left(1+\tfrac{R_1^l}{\lambda_2}\right)
\left(E-\lambda_2+\tfrac{\Lambda_1^l}{2}k\right)
-\tfrac{\Lambda_1^l}{2}\left(1+\tfrac{R_1^l}{\lambda_1}\right)
\left(\tfrac{\Lambda_1^l}{2}+(E+\lambda_2)k\right),
\end{align*}
and
\begin{align*}
\lambda_1=\sqrt{R_1^l(R_1^l+\Lambda_1^l)},\qquad
\lambda_2=\sqrt{R_2^l(R_2^l+\Lambda_2^l)}.
\end{align*}
The derivation is in parallel to the continuous velocity case and its details are left to the Supplemental Material SM2.

\begin{remark}
From the limiting equations in \eqref{eq:p-0-l-two-v-3s}-\eqref{eq:p-0-r-two-v-3s}, summing the $p^{+}$-equation and the $p^{-}$-equation on each subdomain causes the reaction terms on the right-hand side to cancel exactly. Indeed,  we have
\begin{align*}
    &\partial_t(p^{+} + p^{-}) + \partial_x(p^{+}-p^{-}) = 0, \quad x\in\Omega_l\\
&\partial_t(p^{+}+p^{-}) + \partial_x(p^{+}-p^{-}) = 0, \quad x\in\Omega_r.
\end{align*}
 With the interface condition derived by the formal asymptotic analysis, namely
\begin{align} \label{conti:flux}
p^{+}(0^+) - p^{-}(0^+) = p^{+}(0^-) - p^{-}(0^-),
\end{align}
we obtain the flux continuity across $x=0$ for the limiting system. This flux continuity, together with decay at infinity, ensures that the total mass is conserved on the whole space for the limiting system.
\end{remark}

\begin{remark}
Let $A$ be the transition matrix of general $(n+1)$-species such that 
\begin{align} \label{cond:interface}
\begin{pmatrix}
p^{+}(0^+)\\[1mm]
p^{-}(0^+)
\end{pmatrix}
=
A
\begin{pmatrix}
p^{+}(0^-)\\[1mm]
p^{-}(0^-)
\end{pmatrix}
=: \begin{pmatrix}
a_{11} & a_{12}\\[1mm]
a_{21} & a_{22}
\end{pmatrix}
\begin{pmatrix}
p^{+}(0^-)\\[1mm]
p^{-}(0^-)
\end{pmatrix}. 
\end{align}
Then $A$ must have the following structure:
\begin{align} \label{structure:A}
 A=
\begin{pmatrix}
1 & \alpha-1\\[1mm]
0 & \alpha
\end{pmatrix}, 
\qquad \alpha > 0.
\end{align}
The reason that $A$ is upper triangular follows from the derivations in~\eqref{def:A-1} and~\eqref{def:A}: first, $\mathscr{A}_n$ in~\eqref{def:A-1} and $\mathscr{A}$ only rely on the part $v > 0$, which corresponds to $v < 0$ after the change of variable $v \to -v$. Second, due to the continuity of flux across the boundary, equation~\eqref{conti:flux} always holds,
which, given the structure of $A$ in~\eqref{structure:A}, is equivalent to
\begin{align*}
  a_{12}-a_{22} = -1.
\end{align*}
\end{remark}

\section{Numerical validation}
\label{sec:numerical-validation}

In this section, we use particle-based simulations to examine the interface connection condition derived in Section 3. The numerical study addresses three questions. First, for discrete internal-state models in parameter regimes consistent with the scale of the interface-layer analysis, we test whether the internal-state dynamics rapidly produces a time-independent connection law before the population density reaches its steady state. Second, we examine whether the resulting coefficients agree with the theoretical prediction. Third, we repeat the simulations for a continuous internal-state model to determine whether the observed connection law is an artifact of the multi-species discretization.

\subsection{SPECS model and numerical setup}
\label{subsec:numerical-method}

In the SPECS model for \textit{E. coli} chemotaxis~\cite{Jiang2010SPECS, Si2012}, the internal state is the intracellular methylation level $m$, and it evolves according to the following ODE
\begin{equation}
  \frac{\mathrm{d}m}{\mathrm{d}t}
  = f\bigl(m-m_c(x)\bigr),
  \qquad
  f(y)=k_R\left(1-\frac{a(y)}{a_0}\right),
  \label{eq:num-adaptation}
\end{equation}
where the receptor activity is
$a(y)=\frac{1}{1+\exp(-N_r\alpha y)}$.
The equilibrium methylation level given an external signal $s(x)$
is
\[
m_c(x)=M_0+\frac{1}{\alpha}
  \log\left(\frac{1+s(x)/K_I}{1+s(x)/K_A}\right).
\]
We take $M_0=1$, $K_I=18.2\,\mu\mathrm{M}$,
$K_A=3\,\mathrm{mM}$, $\alpha=1.7$, $a_0=0.5$, and $N_r=6$ as for wild-type
\emph{E.~coli} \cite{Jiang2010SPECS}.
The tumbling frequency is
\begin{equation}
  \Lambda\bigl(m-m_c(x)\bigr)
  =z_0+\tau^{-1}\left(\frac{a(m-m_c(x))}{a_0}\right)^H,
  \label{eq:num-tumbling}
\end{equation}
with $z_0=0.14\,\mathrm{s}^{-1}$, $\tau=1.16\,\mathrm{s}$, and
$H=10$.

A perfectly sharp chemical interface is difficult to realize experimentally.
We therefore do not aim for a quantitative comparison with wild-type
\emph{E.~coli}. Single-cell experiments on eukaryotic algae
nevertheless reveal qualitatively similar run-and-tumble motility
\cite{Li2025BiasedLevy} and asymmetric responses at dark--light interfaces
\cite{Diehn1975Photophobic, Lei2025}. However, the light-sensing pathways
of eukaryotic algae are not yet understood well enough to determine
all the parameters required by the model. We therefore use the present parameterization mainly to examine the predicted separation of time scales and the resulting interface law, without attempting quantitative agreement with experiments on a particular microorganism.

The computational domain is $[0,400]~\mu\mathrm{m}$, with
\begin{equation}
  m_c(x)=
  \begin{cases}
    M_r, & x\in[100,300]~\mu\mathrm{m},\\
    M_l, & x\in[0,100)\cup(300,400]~\mu\mathrm{m},
  \end{cases}
  \qquad M_l=1.
  \label{eq:num-signal-profile}
\end{equation}
with periodic boundary conditions imposed. We simulate $2\times10^5$ cells
with velocities $v=\pm v_0$, where
$v_0=16.5/\sqrt{2}~\mu\mathrm{m}/\mathrm{s}$, which is the swim speed of \textit{E. coli} cells. Although the simulations use a swimming speed $v_0\neq 1$ while the
analysis considers $v_0=1$, 
the velocity
magnitude can be absorbed by rescaling the spatial variable, and the resulting
dimensionless interface connection matrix is independent of the specific
value of $v_0$.

For the discrete-$m$ model, the range $[M_l,M_r]$ is divided into $N=n+1$ intervals. Each particle carries a species label $i$ and jumps
between neighboring species with rates inherited from the half-grid values of
$f$. In $\Omega_l$, the transition rate from species $i$ to species $i-1$ is
\begin{equation}
  R_i^l=\frac{|f(m_{i-1/2}-M_l)|}{\Delta m},
  \label{eq:num-transition-rate}
\end{equation}
whereas the tumbling frequency of species $i$ is evaluated at the center of the interval.
For the continuous-$m$ model, the internal state evolves directly according to
equation~\eqref{eq:num-adaptation}. The two particle algorithms are given in the
Appendix. Initially, the cell positions are uniformly
distributed over the spatial domain, and the internal states are uniformly
distributed over different species (for the discrete model) or $[M_l,M_r]$ (for the continuous model). We simulate with a very small time step $\Delta t=0.001$~s.

We evaluate the interface law at $x_0=100~\mu\mathrm{m}$. The one-sided
directional densities $p^{\pm}(t,x_0^{\pm})$ and $p^{\pm}(t,x_0^\mp)$ are obtained by counting particles in the adjacent bins
$[x_0-\Delta x-\Delta,x_0-\Delta]$ and $[x_0+\Delta,x_0+\Delta x+\Delta]$, with
$\Delta x=2~\mu\mathrm{m}$ and $\Delta$ being determined by the width of the interface layers. In the two-velocity setting, the inverse form of
equation~\eqref{cond:interface} is
\begin{equation}
  \begin{pmatrix}
    p^+(t,x_0^-)\\[1mm]p^-(t,x_0^-)
  \end{pmatrix}
  =A^{-1}
  \begin{pmatrix}
    p^+(t,x_0^+)\\[1mm]p^-(t,x_0^+)
  \end{pmatrix},
  \quad
  A^{-1}
  = \begin{pmatrix}
     1 & 1/\alpha -1 \\[3pt]
     0 & 1/\alpha
    \end{pmatrix}
  =:
  \begin{pmatrix}
    1 & -A_{12}\\[1mm]0 & A_{22}
  \end{pmatrix}.
  \label{eq:num-interface-matrix}
\end{equation}
Accordingly, we define
\begin{equation}
  C_1(t)=
  \frac{p^+(t,x_0^+)-p^+(t,x_0^-)}{p^-(t,x_0^+)},
  \qquad
  C_2(t)=\frac{p^-(t,x_0^-)}{p^-(t,x_0^+)},
  \label{eq:num-observables}
\end{equation}
for which equation~\eqref{eq:num-interface-matrix} predicts
\begin{align*}
 C_1(t)\sim A_{12}, 
\quad
 C_2(t)\sim A_{22}, 
\quad
 C_1(t) + C_2 (t) \sim 1. 
\end{align*}
We use the inverse form because
approximating the elements of $A$ would require division by the much smaller
density $p^-(t,x_0^-)$ and would therefore amplify particle-sampling noise.

Since a periodic boundary condition is employed, the system will converge to a steady state at long times, and the densities on each side will converge to constants.
To distinguish the formation of the interface law from the much slower
relaxation of the population densities, we fit both $C_i(t)$ ($i=1,2$) and the four
one-sided densities with
\begin{equation}
  E(t)=c+a\exp(\lambda t),\qquad \lambda<0.
  \label{eq:num-exponential-fit}
\end{equation}
The magnitude
$|\lambda|$ measures the corresponding relaxation rate.  

\subsection{Dynamics of the interface problem}
\label{subsec:asymptotic-validation}

\begin{figure}[!htp]
\centering
  \begin{subfigure}{0.36\textwidth}
    \centering
    \includegraphics[width=\textwidth]{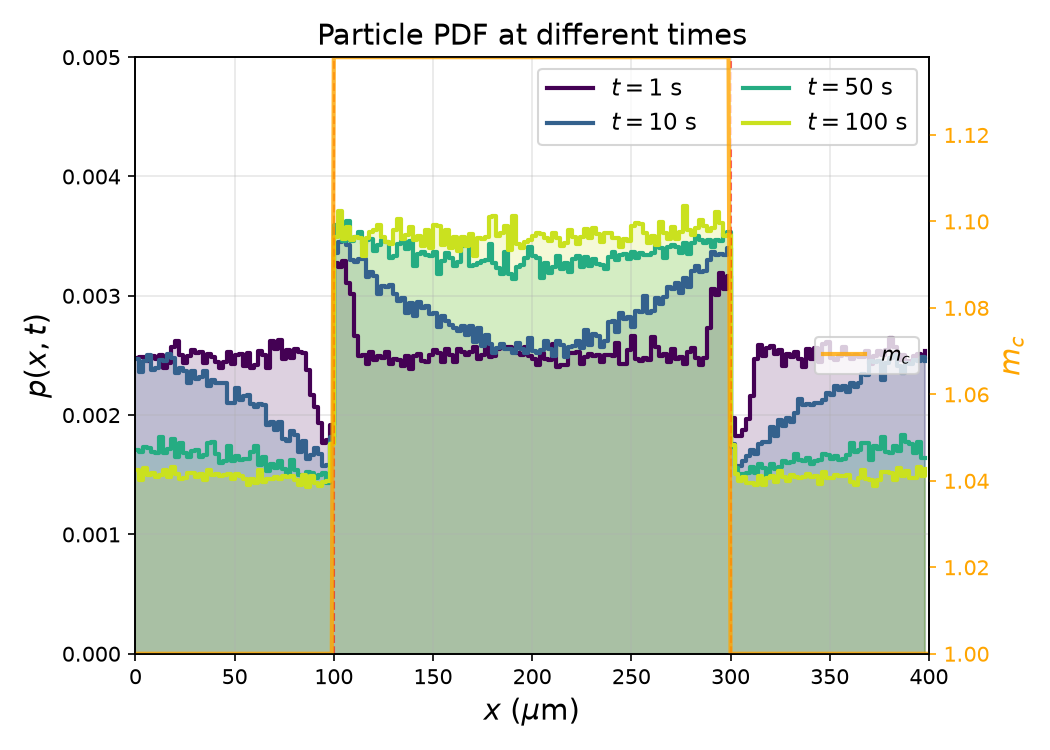}
    \caption{}
  \end{subfigure}
     \begin{subfigure}{0.3\textwidth}
    \centering
    \includegraphics[width=\textwidth]{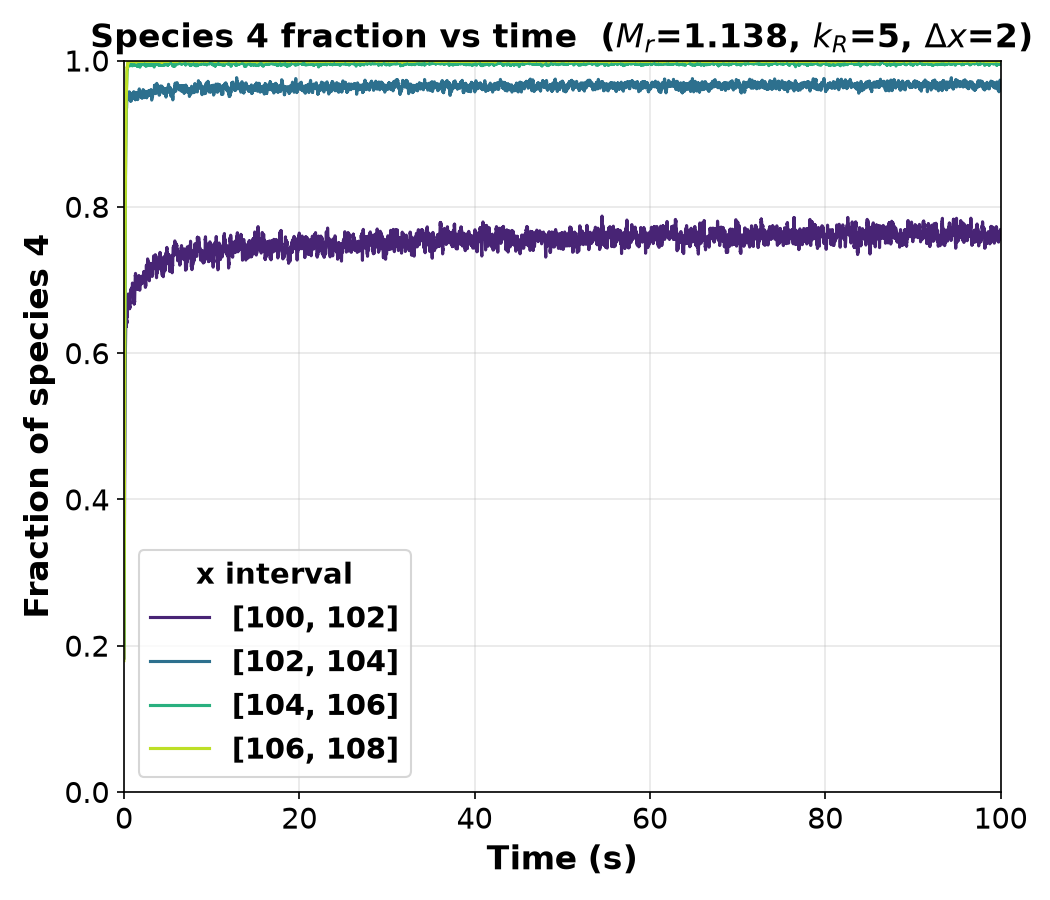}
    \caption{}
  \end{subfigure}
    \begin{subfigure}{0.30\textwidth}
    \centering
    \includegraphics[width=\textwidth]{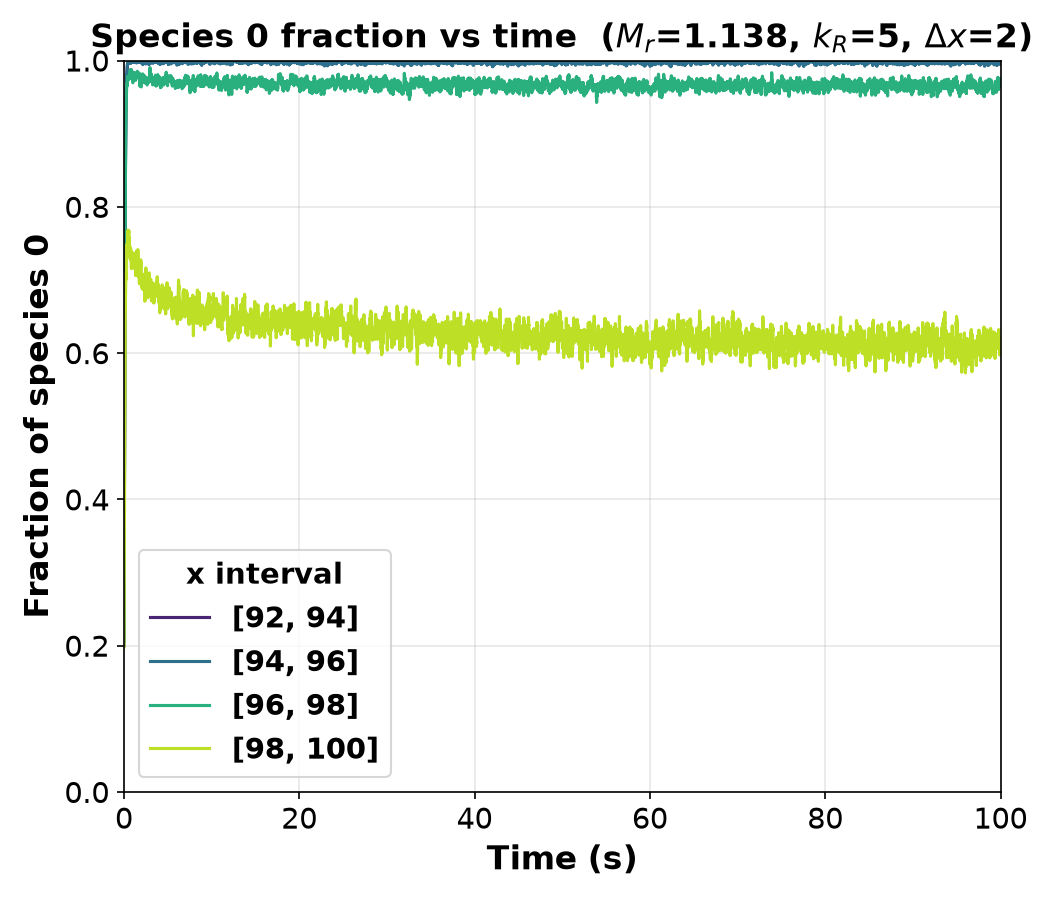}
    \caption{}
  \end{subfigure}
  \begin{subfigure}{0.30\textwidth}
    \centering
    \includegraphics[width=\textwidth]{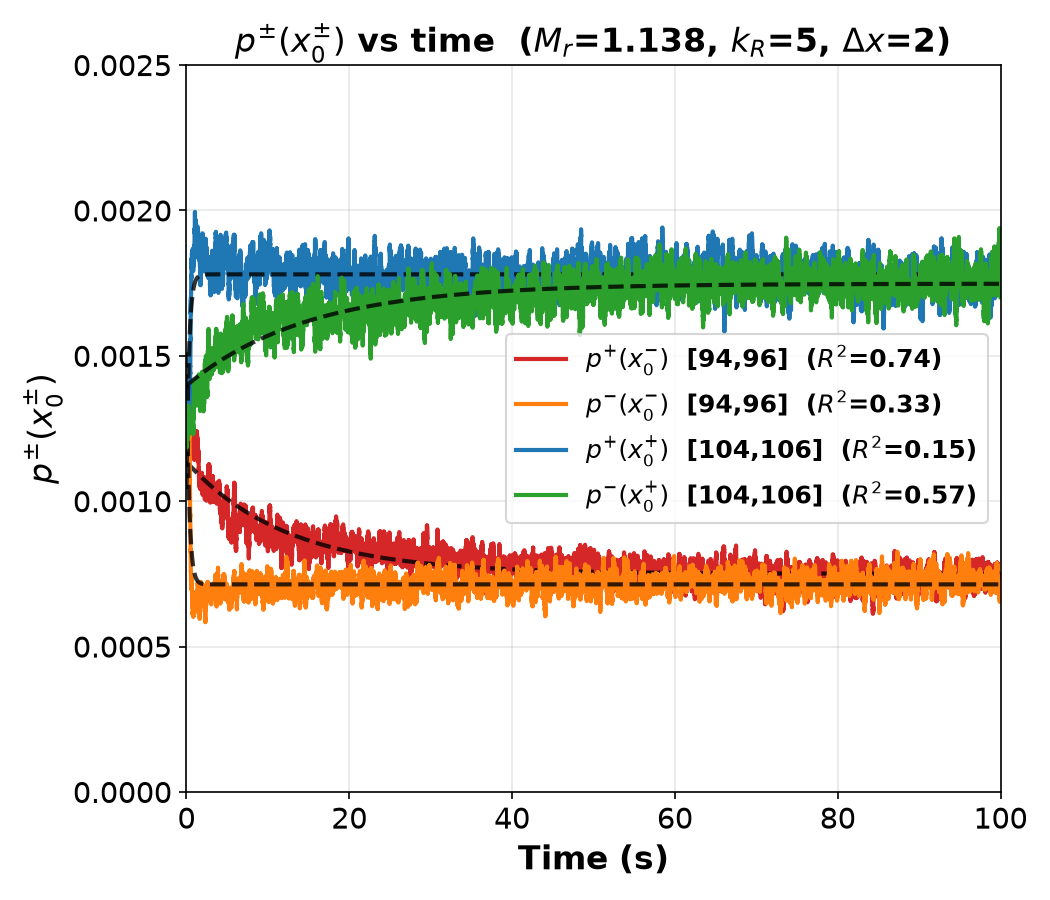}
    \caption{}
  \end{subfigure}
   \begin{subfigure}{0.32\textwidth}
    \centering
    \includegraphics[width=\textwidth]{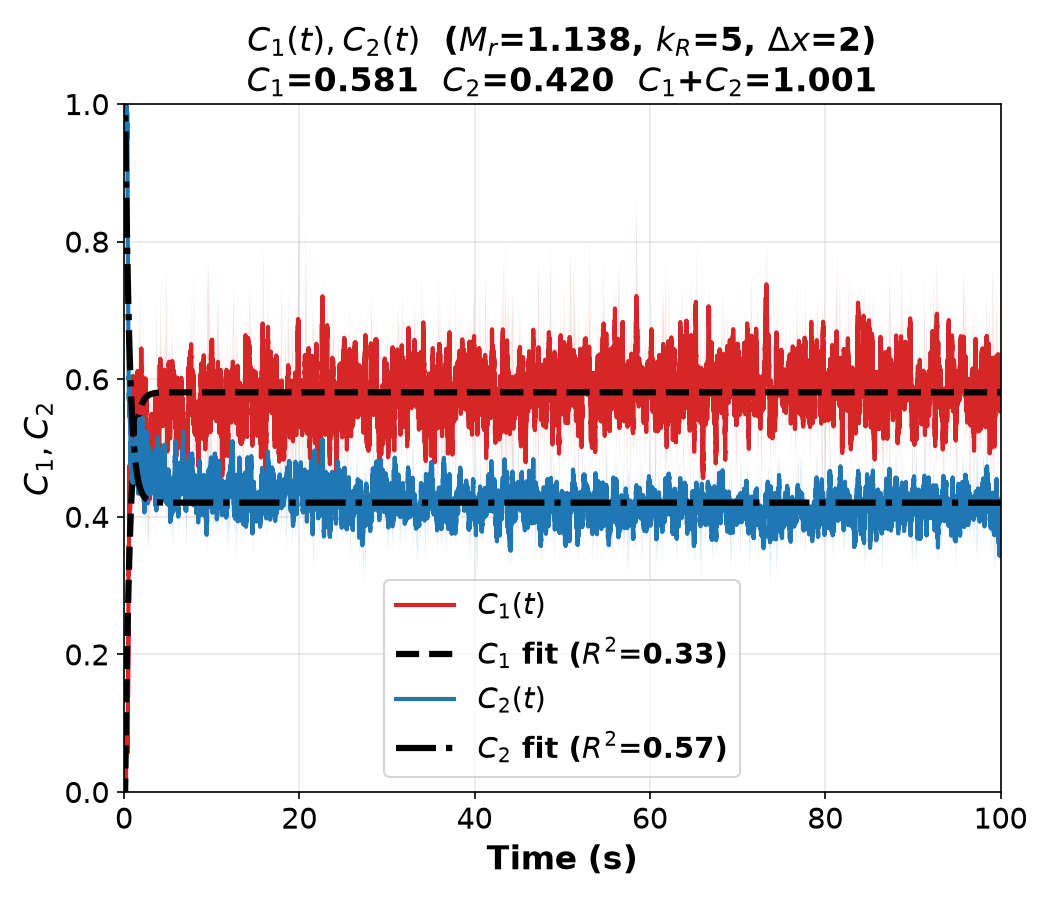}
    \caption{}
  \end{subfigure}
  \caption{Interface dynamics for the five-species model with
$k_R=5~\mathrm{s}^{-1}$ and $M_r=1.138$.
(a) Spatial density profiles at $t=1$, $10$, $50$, and $100~\mathrm{s}$; (b) Fractions of species 4 near the interface at $M_l=1$; (c) Fractions of species 0 near the interface at $M_r=1.138$; 
(d) One-sided directional densities near $x_0=100~\mu\mathrm{m}$:
$p^+(t,x_0^+)$, $p^-(t,x_0^+)$ denote the particle densities inside $[104\mu m,106\mu m]$ and $p^+(t,x_0^-)$, 
$p^-(t,x_0^-)$ are for the interval $[94\mu m,96\mu m]$ . The densities are normalized by the total particle
number divided by the bin width.
(e) $C_1(t)$ and $C_2(t)$, together with their fitted constant values.
Although the density profiles $p^\pm$ continue to evolve,
the interface ratios rapidly approach constant values. }
  \label{fig:num-original-kr05}
\end{figure}

We consider the following five-species case:
\begin{align}
  (k_R,M_l,M_r)
  &=
  (5~\mathrm{s}^{-1},1,1.138).
  \label{eq:num-two-asymptotic-cases}
\end{align}
The corresponding transition and tumbling rates are
$$
\begin{aligned}
&(R_0^l,R_1^l,R_2^l,R_3^l,R_4^l)=(0,12.72,37.39,59.93, 79.45),\\
&(\Lambda_0^l,\Lambda_1^l,\Lambda_2^l,\Lambda_3^l,\Lambda_4^l)=(1.002,4.433,16.00,44.72,99.02).
\end{aligned}
$$
Thus, $R_n^l$ and $\Lambda_n^l$ are comparable as required by the asymptotic scaling.
The numerical results obtained with $\Delta t=0.001~\mathrm{s}$ are
shown in Fig.~\ref{fig:num-original-kr05}.

Simulations with this parameter set assess the biological relevance of the asymptotic framework in \eqref{sys:n-species-0}--\eqref{sys:n-species-N}. The observed interface dynamics are qualitatively consistent with the predictions of the asymptotic analysis.
The spatial density profiles in Fig.~\ref{fig:num-original-kr05}(a)
show a continuous redistribution of the population toward the
high-signal region $\Omega_r$. Even after $t=50~\mathrm{s}$, the density
in $\Omega_r$ continues to increase, while in $\Omega_l$ it continues
to decrease.

Figure~\ref{fig:num-original-kr05}(b)(c) 
show that the internal-state
distribution relaxes rapidly and becomes concentrated near species
$p_0$ inside $\Omega_l$ and species $p_4$ inside $\Omega_r$. 
 The interface layer extends approximately $4\,\mu\mathrm{m}$ into each
adjacent region. Outside these layers, the locally adapted species
dominates: $p_0$, associated with $m_0 = 1.000$, in the low-signal regions,
and $p_4$, associated with $m_4=1.138$,  in the high-signal region.
Within the interface layers, however, the contributions from the
other species remain appreciable.
These observations support an RTKM description in the bulk regions $[0,96]\cup[304,400]$ and $[104,296]$, with all spatial coordinates
expressed in micrometers. The bulk solutions are coupled through
effective connection conditions that account for the dynamics within
the interface layers. To estimate the interface coefficients
at $x=100\,\mu\mathrm{m}$, we use the simulated densities
averaged over $[94,96]$ and $[104,106]$, which are immediately outside the layer
on its left and right sides, respectively. These sampling windows
approximate the outer densities entering the interface condition.
An analogous procedure applies at $x=300\,\mu\mathrm{m}$.

For the interface condition, the separation of time scales is
demonstrated more directly by comparing
Figs.~\ref{fig:num-original-kr05}(d) and
\ref{fig:num-original-kr05}(e).
The four one-sided densities in Figure 1(d) remain visibly time-dependent. In particular, particles moving away from the interface quickly reach a steady state and the rate is at 
$|\lambda_{p^\pm(x_0^\pm)}|
  =
  3.1~\mathrm{s}^{-1}$
whereas particles moving toward the interface relax slowly with
$|\lambda_{p^\pm(x_0^\mp)}|
  =
  0.070~\mathrm{s}^{-1}$.
By contrast, the combinations $C_1(t)$ and $C_2(t)$ shown in Figure 1(e)
rapidly approach constants, with fitted relaxation rates 
$  |\lambda_{C_i}|
  =
  1.9~\mathrm{s}^{-1}$.
Thus, the interface coefficients relax more than one order of magnitude
faster than the dominant one-sided densities. Consequently, the relation
between the incoming and outgoing densities becomes effectively
time-independent long before the population density near the interface
reaches its steady state. This separation of time scales creates a
substantial intermediate-time window in which the population dynamics can
already be described by the RTKM supplemented by the time-independent
interface condition \eqref{eq:num-interface-matrix}, even though the population itself remains
far from steady state.

\subsection{Influence of intermediate-state tumbling on the interface coefficients}
\label{subsec:tumbling-controls}

Having observed the rapid relaxation of the interface ratios, we now
examine their dependence on the tumbling rates for
$k_R=5~\mathrm{s}^{-1}$, $M_l=1$, and $M_r=1.138$.
We denote the fitted plateaus of $C_1(t)$ and $C_2(t)$ by $c_1$ and
$c_2$, respectively. For these parameters,
\[
\frac{\Lambda_1^l}{R_1^l}\approx0.349,\qquad
\frac{\Lambda_2^l}{R_2^l}\approx0.428,\qquad
\frac{\Lambda_3^l}{R_3^l}\approx0.746.
\]
These ratios indicate appreciable reorientation during adaptation
through the intermediate states.
Table~\ref{tab:interface-coefficients} compares the baseline model with
two controls: the first sets $\Lambda_1^l=\Lambda_2^l=\Lambda_3^l=0$,
and the second additionally sets $\Lambda_0^l=0$. In the particle
simulations and the full stationary model, all right-side rates remain
unchanged. The three numerical columns distinguish the fitted particle
results, the interface-layer prediction, and the full stationary
prediction evaluated using the simulation sampling windows.

\paragraph{Interface-layer prediction}
The second numerical column is obtained from the two-velocity version
of \eqref{eq:BL-L} and \eqref{eq:BL-R} with five species.
After expressing both half-space problems in a common coordinate
$\xi\in\mathbb R$, oriented in the physical $x$ direction, and
absorbing the velocity magnitude into $\xi$, the ten directional
densities $q_i^\pm(\xi)$ satisfy
\begin{subequations}\label{eq:num-layer-system}
\begin{align}
\pm\frac{\mathrm d q_i^\pm}{\mathrm d\xi}
&=\frac{\ell_i^l}{2}(q_i^\mp-q_i^\pm)
  +R_{i+1}^lq_{i+1}^\pm-R_i^lq_i^\pm,
&&\xi<0,\label{eq:num-layer-left}\\
\pm\frac{\mathrm d q_i^\pm}{\mathrm d\xi}
&=R_{i-1}^rq_{i-1}^\pm-R_i^rq_i^\pm,
&&\xi>0,\label{eq:num-layer-right}
\end{align}
\end{subequations}
for $i=0,\ldots,4$. Here $R_0^l=R_4^r=0$, terms involving
species outside $\{0,\ldots,4\}$ are zero, and
\[
\ell_0^l=0,\qquad
\ell_i^l=\Lambda_i^l\quad(i=1,\ldots,4),
\]
with the stated control modifications applied to $\ell_1^l,\ell_2^l,
\ell_3^l$. In particular, the layer system omits adapted-state
tumbling on the left and all tumbling on the right, as prescribed
by the leading-order half-space equations. These omissions distinguish
it from the full stationary system.

We impose continuity of every directional species at $\xi=0$ and
boundedness as $\xi\to\pm\infty$. Computing the following far-field total density 
\[
P_l^\pm=\lim_{\xi\to-\infty}\sum_{i=0}^4q_i^\pm(\xi),
\qquad
P_r^\pm=\lim_{\xi\to+\infty}\sum_{i=0}^4q_i^\pm(\xi)
\]is equivalent to the derivation of the interface connection matrix in Section \ref{subsec:three-species}.
The baseline layer system gives
$(C_1,C_2)
=(0.605866,0.394134)$.
For the two settings with $\Lambda_i = 0$, it gives $(0.332776,0.667224)$.
The agreement of the two predictions follows because
$\Lambda_0^l$ is already absent from the layer system.

\paragraph{Full stationary prediction} 
The third numerical column is calculated independently from the
stationary five-species kinetic equations in the original periodic
domain $[0,400]~\mu\mathrm m$. Writing the stationary directional
densities as $u_i^\pm(x)$, we solve
\begin{subequations}\label{eq:num-full-stationary}
\begin{align}
\pm v_0\frac{\mathrm d u_i^\pm}{\mathrm d x}
&=\frac{\Lambda_i^l}{2}(u_i^\mp-u_i^\pm)
 +R_{i+1}^lu_{i+1}^\pm-R_i^lu_i^\pm,
&&x\in\Omega_l,\\
\pm v_0\frac{\mathrm d u_i^\pm}{\mathrm d x}
&=\frac{\Lambda_i^r}{2}(u_i^\mp-u_i^\pm)
 +R_{i-1}^ru_{i-1}^\pm-R_i^ru_i^\pm,
&&x\in\Omega_r,
\end{align}
\end{subequations}
with the same endpoint conventions, periodicity, and continuity of
all species at $x=100$ and $300~\mu\mathrm m$. All tumbling rates
are retained except those explicitly suppressed in each control.

For $h=2~\mu\mathrm m$, we evaluate
\[
\overline P_l^\pm=\frac1h\int_{94\,\mu\mathrm m}^{96\,\mu\mathrm m}
 \sum_{i=0}^4u_i^\pm(x)\,\mathrm d x,
\qquad
\overline P_r^\pm=\frac1h\int_{104\,\mu\mathrm m}^{106\,\mu\mathrm m}
 \sum_{i=0}^4u_i^\pm(x)\,\mathrm d x,
\]
and get 
$$
C_1^{\mathrm{S}}=
\frac{\overline P_r^+-\overline P_l^+}{\overline P_r^-},
\qquad
C_2^{\mathrm{S}}=\frac{\overline P_l^-}{\overline P_r^-}.
$$
These are stationary predictions for the same bin-based observables
used in the particle simulations. They are distinct from the
far-field coefficients of the leading-order interface-layer problem.

\begin{table}[htbp]
\centering
\small
\setlength{\tabcolsep}{4pt}
\begin{tabular}{lccc}
\hline
Model or control
& \shortstack{Particle\\simulations}
& \shortstack{Interface layer\\far-field values}
& \shortstack{Full stationary ODE\\bin averages}\\
\hline
Baseline
& $(0.581,\,0.420)$
& $(0.605866,\,0.394134)$
& $(0.595284,\,0.404716)$\\
$\Lambda_1^l=\Lambda_2^l=\Lambda_3^l=0$
& $(0.290,\,0.711)$
& $(0.332776,\,0.667224)$
& $(0.314503,\,0.685497)$\\
$\Lambda_0^l=\Lambda_1^l=\Lambda_2^l=\Lambda_3^l=0$
& $(0.344,\,0.659)$
& $(0.332776,\,0.667224)$
& $(0.356253,\,0.643747)$\\
\hline
\end{tabular}
\caption{Interface ratios $(C_1,C_2)$ for $k_R=5~\mathrm{s}^{-1}$,
$M_l=1$, and $M_r=1.138$. The three numerical columns report,
respectively, fitted particle-simulation plateaus, coefficients
obtained from the bounded whole-line interface-layer solution
\eqref{eq:num-layer-system}, and stationary bin-based ratios from
\eqref{eq:num-full-stationary}.}
\label{tab:interface-coefficients}
\end{table}

The interface-layer calculation captures the substantial decrease
in $C_1$ and increase in $C_2$ caused by suppressing intermediate-state
tumbling. However, it predicts identical coefficients for the two
controls, whereas the particle simulations and the full stationary
system both show a change when $\Lambda_0^l$ is additionally removed.
This difference identifies a finite-parameter contribution from a
process omitted in the leading-order layer equations. Together, the comparisons show that the
interface-layer analysis provides a time-independent scattering law
whose coefficients depend on the left-side tumbling processes
retained at leading order, while processes omitted in that limit
can modify the finite-parameter stationary ratios.

\subsection{Continuous internal state}
\label{subsec:continuous-m}

Finally, we replace the discrete methylation states with a continuous
internal variable $m$, evolving according to
\eqref{eq:num-adaptation}. We retain $M_l=1$ and consider
$M_r=1.138$ and $1.076$, with all other parameters and the
density-extraction procedure unchanged.

\begin{figure}[t]
  \centering
  \begin{subfigure}{0.36\textwidth}
    \centering
    \includegraphics[width=\textwidth]{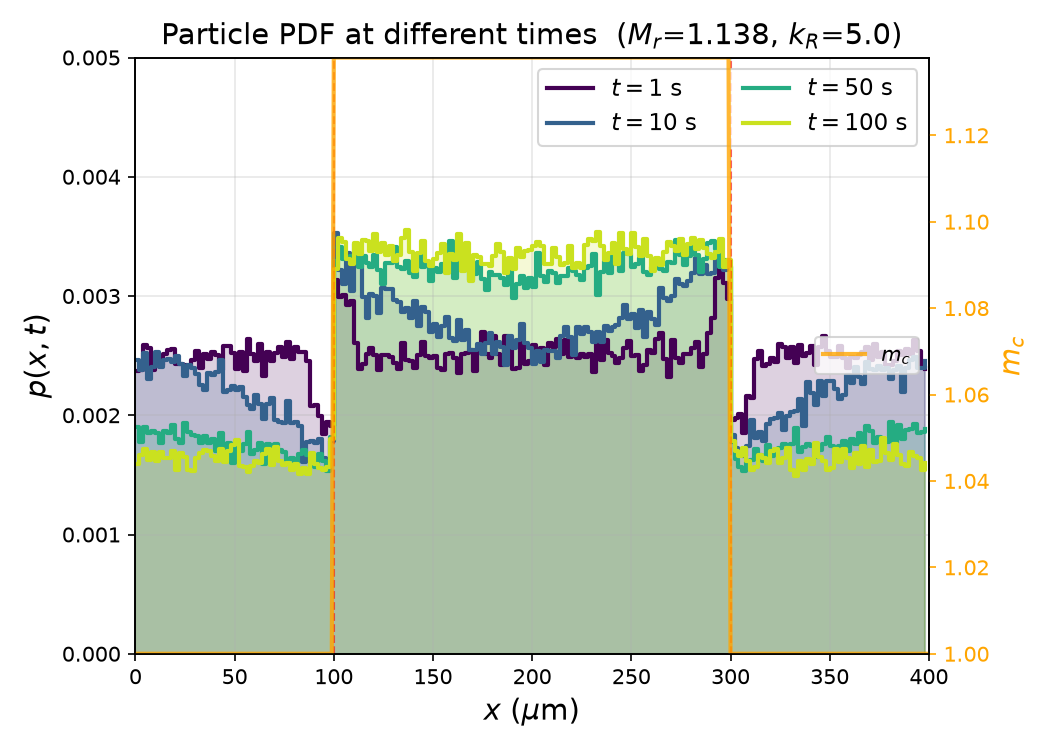}
    \caption{}
  \end{subfigure}
  \begin{subfigure}{0.30\textwidth}
    \centering
    \includegraphics[width=\textwidth]{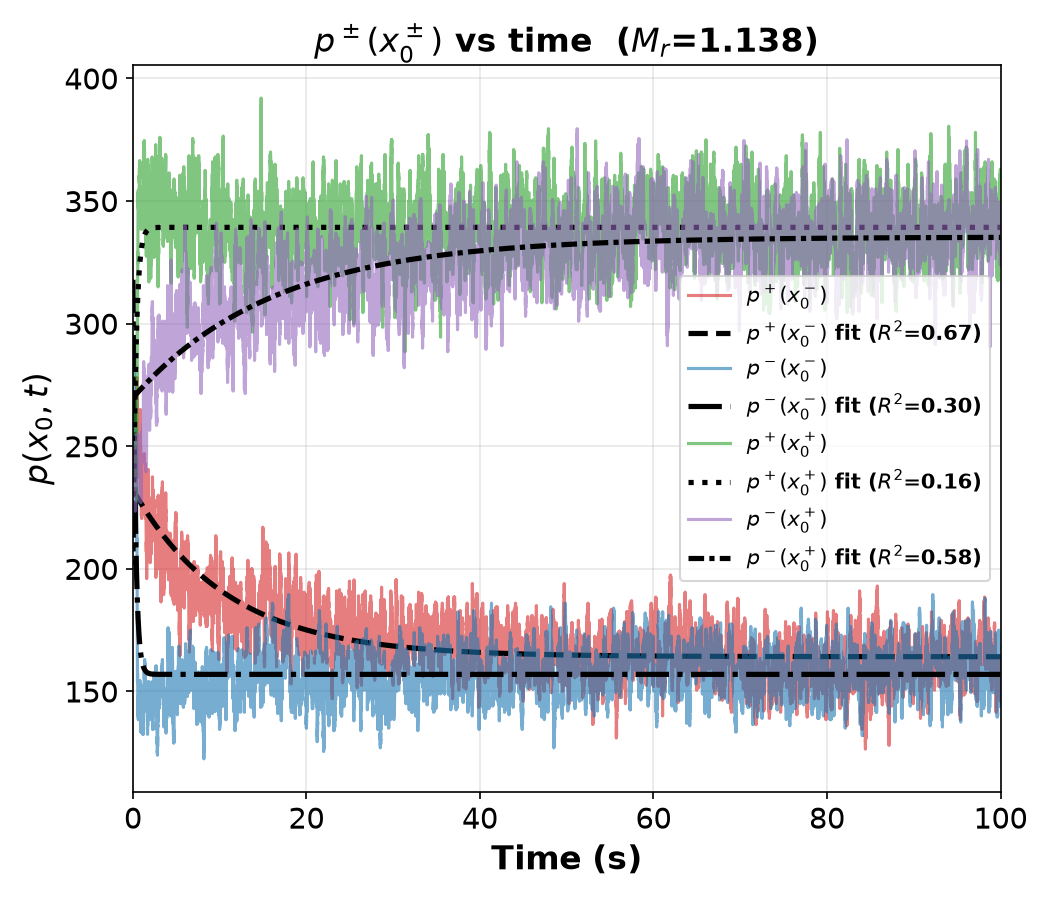}
    \caption{}
  \end{subfigure}
  \begin{subfigure}{0.30\textwidth}
    \centering
    \includegraphics[width=\textwidth]{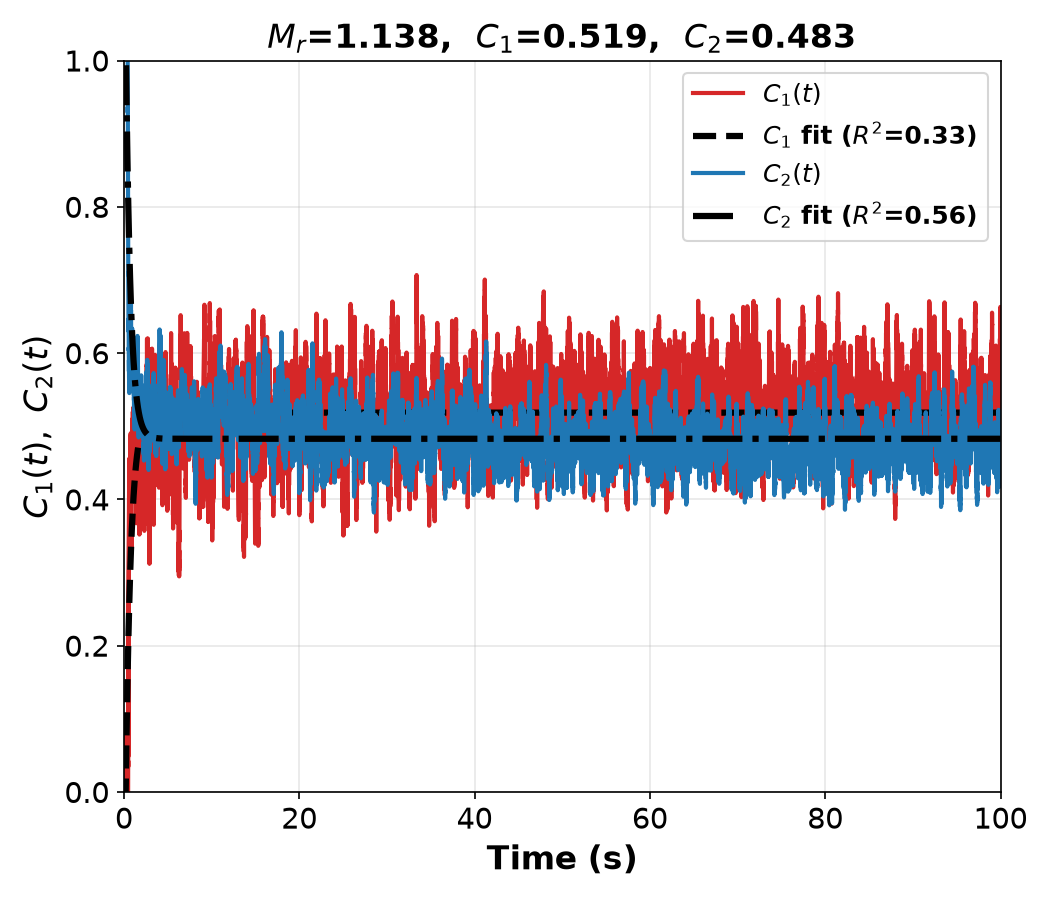}
    \caption{}
  \end{subfigure}
  
  \begin{subfigure}{0.36\textwidth}
    \centering
    \includegraphics[width=\textwidth]{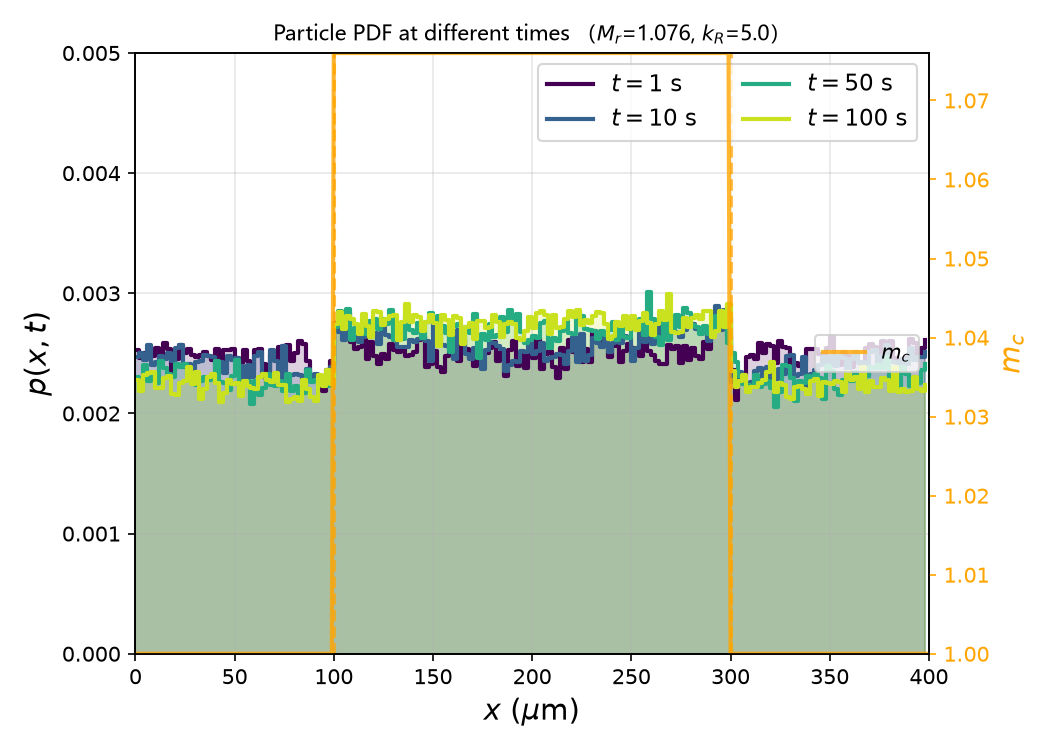}
    \caption{}
  \end{subfigure}
  \begin{subfigure}{0.30\textwidth}
    \centering
    \includegraphics[width=\textwidth]{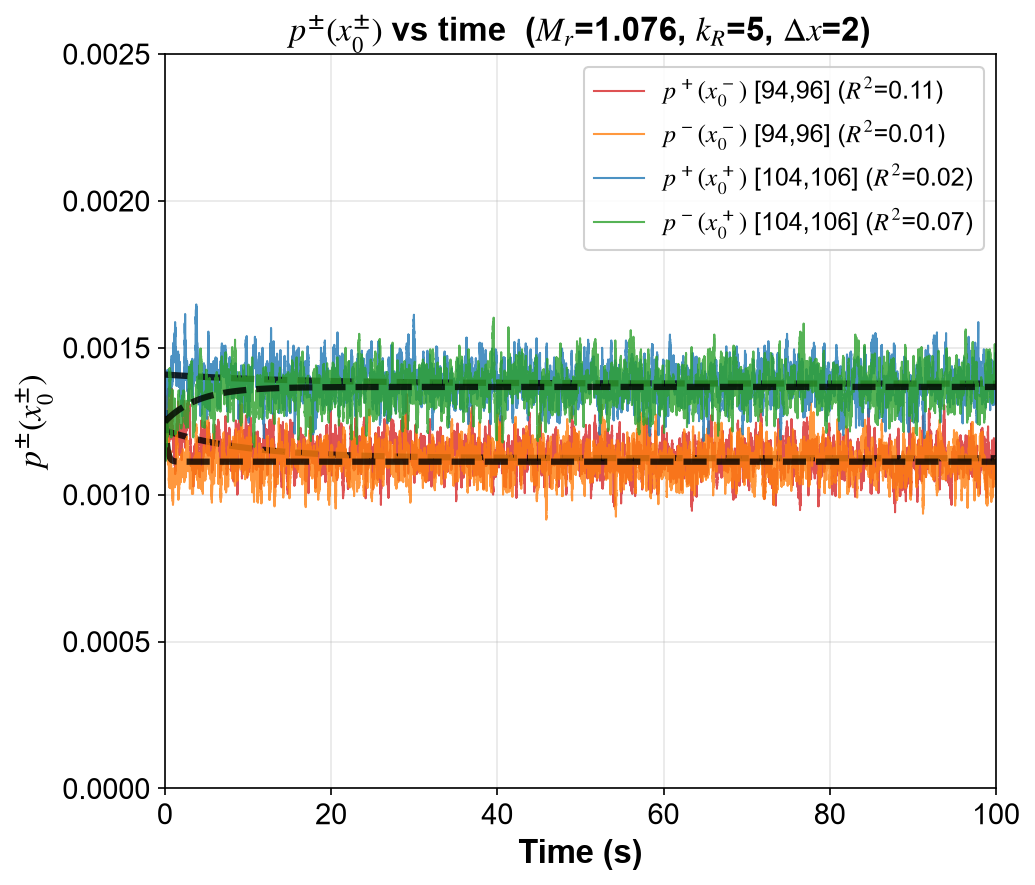}
    \caption{}
  \end{subfigure}
  \begin{subfigure}{0.30\textwidth}
    \centering
    \includegraphics[width=\textwidth]{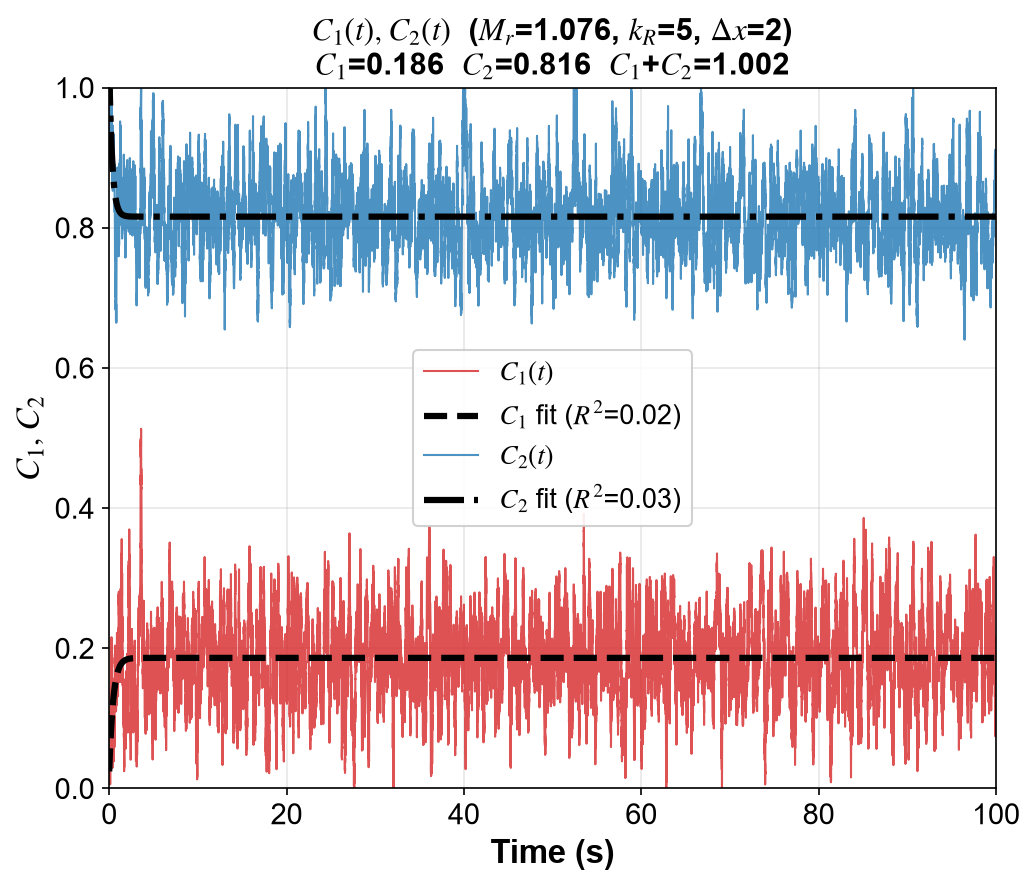}
    \caption{}
  \end{subfigure}
  \caption{Influence of $M_r$ on the interface dynamics in the pathway-based kinetic
  model with $M_l=1$ $M_r=1.138$ (top row) and $M_l=1$ $M_r=1.076$ (bottom row).
  (a) (d): Spatial density profiles at $t=1$, $10$, $50$, and $100~\mathrm{s}$;
(b) (e) One-sided directional densities near $x_0=100~\mu\mathrm{m}$:
$p^+(t,x_0^+)$, $p^-(t,x_0^+)$ denote the particle densities inside $[104\mu m,106\mu m]$ and $p^+(t,x_0^-)$, 
$p^-(t,x_0^-)$ are for the interval $[94\mu m,96\mu m]$ . The densities are normalized by the total particle
number divided by the bin width.
(c) (f) $C_1(t)$ and $C_2(t)$, together with their fitted constant values.}
  \label{fig:num-original-continuous}
\end{figure}

In Figure~\ref{fig:num-original-continuous}, the left column shows the redistribution of cells toward the central
high-signal region $\Omega_r$, with stronger accumulation for larger
$M_r$. The middle column shows the time evolution of the four
one-sided directional densities, measured in the sampling intervals
$[94,96]\,\mu\mathrm{m}$ and $[104,106]\,\mu\mathrm{m}$.
The right column shows $C_1(t)$ and $C_2(t)$, both of which approach
nearly constant plateaus after a short transient. This observation indicates that the rapid establishment of an approximately
time-independent interface relation is not an artifact of internal-state discretization.

For $M_r=1.138$, the fitted plateau values are also close to those
obtained from the five-species model, suggesting that this discrete
model reasonably approximates the interface ratios of the continuous
PBKM for the parameters considered. Deriving an asymptotic connection
law directly from the continuous PBKM remains an interesting challenge.
In particular, the interplay between rapid internal-state adaptation
and the strong dependence of the tumbling rate on $m$ complicates
the identification of suitable asymptotic scales.

\section{Conclusion and discussion}\label{sec:conclusion}
In this work, we derive an interface scattering law for run-and-tumble 
microorganisms whose responses to abrupt environmental changes are mediated by 
intracellular biochemical adaptation. Starting from a PBKM, asymptotic analysis yields RTKMs in the two 
homogeneous regions, coupled by an effective scattering condition at the 
dark--light interface. In slab geometry, this condition takes the form of an integral operator 
in velocity space, while for the two-velocity model the operator reduces to an explicit transition
matrix. Unlike conventional optical interface laws, the outgoing 
distribution at one velocity can depend on the entire incoming velocity 
distribution because of adaptation and tumbling within the interface layer.

Numerical simulations of the two-flux PBKM based on biologically relevant parameters show that a time-independent interface law emerges before the population reaches 
steady state. The same structure is also observed for a continuous internal state. This suggests that 
effective interface scattering 
may be generic in microorganisms that respond to environmental signals through intracellular adaptation.

Our analysis provides theoretical support for the idea that interface scattering laws may be tuned by modifying receptors, adaptation 
enzymes, or feedback loops. Microbial transport across environmental interfaces may therefore be controlled through modifications of intracellular signaling pathways.

\newpage
\bibliographystyle{siamplain}
\bibliography{references}
\end{document}


\maketitle

\section{An elementary technical lemma}

\begin{lemma} \label{lem:expo-decay}
Suppose $f$ is a bounded solution to the ODE
\begin{align} \label{ODE}
  v \p_x f + A f(x, v) = Bg(x, v), 
\qquad
  \text{$f(0, v)$ bounded,}
\end{align}
where $x \in [0, \infty)$, $v \in (-1, 1)$, $A, B$ are constants, $A \neq 0$, and $|g (x, v)| \leq c_1 e^{-\alpha x}$ for some constants $c_1, \alpha > 0$. Then for any $v \in (-1, 1)$, 
\begin{align*}
  |f(x, v)| \to 0, 
\qquad
  \text{as $x \to \infty$.}
\end{align*}
\end{lemma} 
\begin{proof}
For $v \neq 0$, denote 
\begin{align*}
  A_1 =  \frac{A}{v}, 
\qquad
  B_1 = \frac{B}{v}, 
\end{align*}
where both $A_1, B_1$ are independent of $x$. Then $f$ solves
\begin{align} \label{ODE:1}
  \p_x f + A_1 f = B_1 g. 
\end{align}
Directly solving~\eqref{ODE:1} gives
\begin{align*}
  f(x, v)
= f(0, v)e^{-A_1 x}
  + B_1 \int_0^x e^{A_1 (s-x)} g(s, v) \d s. 
\end{align*}
If $A_1 > 0$, then 
\begin{align*}
  |f(x, v)|
&\leq
 (\sup_v |f(0, v)|) \, e^{-A_1 x}
 + |B_1| c_1 \int_0^x e^{-A_1 (x-s)} e^{-\alpha s} \d s
\\
&\leq
(\sup_v |f(0, v)|) \, e^{-A_1 x}
 + |B_1| c_1 \int_0^{x/2} e^{-A_1 (x-s)} e^{-\alpha s} \d s
 \\
&\quad +|B_1| c_1 \int_{x/2}^x e^{-A_1 (x-s)} e^{-\alpha s} \d s
\\
&\leq
  (\sup_v |f(0, v)|) \, e^{-A_1 x}
  + C e^{-\frac{A_1}{2} x}
  + C e^{-\frac{\alpha}{2} x}
\to 0
\qquad \text{as $x \to \infty$.}
\end{align*}
If $A_1 < 0$, then 
\begin{align*}
  f(x, v)
= e^{-A_1 x}
  \vpran{f(0, v) + B_1 \int_0^x e^{A_1 s} g(s, v) \d s}.
\end{align*}
Note that since $g$ decays exponentially, the integral converges. Therefore, there exists $E_0$ such that 
\begin{align*}
  f(x, v) \sim E_0 e^{-A_1 x}
\qquad
 \text{for $x \gg 1$.} 
\end{align*}
Given the boundedness of $f$, one must have $E_0 = 0$. Therefore, we can reformulate $f(x, v)$ as
\begin{align*}
  f(x, v)
= - e^{-A_1 x}
  \int_x^\infty e^{A_1 s} g(s, v) \d s,
\end{align*}
which gives
\begin{align*}
  |f(x, v)|
&\leq
  c_1 \int_x^\infty e^{A_1 (s -x)} e^{-\alpha s} \d s
\\
&\leq
 c_1 \int_x^\infty e^{-\alpha s} \d s
 = \frac{c_1}{\alpha} e^{-\alpha x} \to 0,
\qquad
 \text{as $x \to \infty$.}
\end{align*}
Hence, regardless of the sign of $A$, we always have $f(x, v) \to 0$ as $x \to \infty$.
\end{proof}

\section{The computation details of the transition matrix $A_2$}

The interface operator $\mathscr{A}$ in this setting will be a $2\times 2$ transition matrix, which depends on the solution operator  of the $6\times 6$ ODE system. 

\medskip

\noindent\textbf{Left domain $\Omega_l$.} In the left domain $\Omega_l$, equation (3.7) becomes 
\begin{align*}
\Omega_l: \quad
\begin{cases}
\partial_y\widetilde p_0^+= R_1^l\widetilde p_1^+,
\\
\partial_y\widetilde p_0^-
= -R_1^l\widetilde p_1^-,
\\
\partial_y\widetilde p_1^{+} 
+ (R_1^l + \tfrac{\Lambda_1^l}{2})\widetilde p_1^{+} - \tfrac{\Lambda_1^l}{2}\widetilde p_1^{-} 
= R_2^l\widetilde p_2^{+},
\\
\partial_y\widetilde p_1^{-} 
- (R_1^l + \tfrac{\Lambda_1^l}{2})\widetilde p_1^{-} + \tfrac{\Lambda_1^l}{2}\widetilde p_1^{+} 
= -  R_2^l\widetilde p_2^{-},
\\
 \partial_y \widetilde p^{+}_{2}
 + R_2^l \widetilde p^{+}_{2}
  = \tfrac12\Lambda_{2}^l (\widetilde p^{-}_{2} - \widetilde p^{+}_{2}),
\\
 \partial_y \widetilde p^{-}_{2}
 - R_2^l \widetilde p^{-}_{2}
 = \tfrac12\Lambda_{2}^l (\widetilde p^{-}_{2} - \widetilde p^{+}_{2}),
\end{cases}
\end{align*}
which can be solved by solving subsystems for $\widetilde p_2^\pm, \widetilde p_1^\pm$, and $\widetilde p_0^\pm$. First, the equations for $\widetilde p_2^\pm$ form a closed system: 
\begin{align*}
\begin{cases}
 \partial_y \widetilde p^{+}_{2}
 + R_2^l \widetilde p^{+}_{2}
 = \tfrac12\Lambda_{2}^l (\widetilde p^{-}_{2} - \widetilde p^{+}_{2}),
\\
 \partial_y \widetilde p^{-}_{2}
 - R_2^l \widetilde p^{-}_{2}
 = \tfrac12\Lambda_{2}^l (\widetilde p^{-}_{2} - \widetilde p^{+}_{2}).
\end{cases}
\end{align*}
The eigenvalues of the $2\times 2$ matrix of the reaction terms are $\pm\lambda_2$ with
\begin{align*}
\lambda_2=\sqrt{R_2^l(R_2^l+\Lambda_2^l)}.
\end{align*}
By the boundedness of the solution as $y\to\infty$, we have
\begin{align} \label{soln:tilde-p-2}
\widetilde p_2^+(t,y)=\widetilde C_2(t)e^{-\lambda_2 y},
\qquad
\widetilde p_2^-(t,y)
=\frac{\Lambda_2^l+2R_2^l-2\lambda_2}{\Lambda_2^l}
\widetilde C_2(t)e^{-\lambda_2 y}.
\end{align}
Substituting this into the $\widetilde p_1^\pm$-equations
\begin{align*}
\begin{cases}
\partial_y\widetilde p_1^{+} = - (R_1^l + \tfrac{\Lambda_1^l}{2})\widetilde p_1^{+} + \tfrac{\Lambda_1^l}{2}\widetilde p_1^{-} 
+ R_2^l\widetilde p_2^{+},
\\
\partial_y\widetilde p_1^{-} 
= (R_1^l + \tfrac{\Lambda_1^l}{2})\widetilde p_1^{-} - \tfrac{\Lambda_1^l}{2}\widetilde p_1^{+} -  R_2^l\widetilde p_2^{-},
\end{cases}
\end{align*}
whose eigenvalues of the reaction matrix are $\pm\lambda_1$ with
\begin{align*}
\lambda_1=\sqrt{R_1^l(R_1^l+\Lambda_1^l)}.
\end{align*}
Here we need to distinguish two cases: $\lambda_1 = \lambda_2$ and $\lambda_1 \neq \lambda_2$. In the case where $\lambda_1 \neq \lambda_2$, by the boundedness of $\widetilde p_2^\pm$ and $\widetilde p_1^\pm$, we have
\begin{align} 
\widetilde p_1^+(t,y)
=&\tfrac{\Lambda_1^l}{2}
\widetilde C_1(t)e^{-\lambda_1 y} \no
\\
&+\tfrac{R_2^l\widetilde C_2(t)}
{(R_1^l)^2+ R_1^l\Lambda_1^l -\lambda_2^2}
\left[
R_1^l+\tfrac{\Lambda_1^l}{2}+\lambda_2
+\tfrac{\Lambda_1^l}{2}
\tfrac{\Lambda_2^l+2R_2^l-2\lambda_2}{\Lambda_2^l}
\right]e^{-\lambda_2 y},
\label{soln:tilde-p1-p}
\\[3mm]
\widetilde p_1^-(t,y)
=&\left(R_1^l+\tfrac{\Lambda_1^l}{2} - \lambda_1\right)
\widetilde C_1(t)e^{-\lambda_1 y} \no
\\
&+\tfrac{R_2^l\widetilde C_2(t)}
{(R_1^l)^2+ R_1^l\Lambda_1^l-\lambda_2^2}
\left[
\tfrac{\Lambda_1^l}{2}
+\left(R_1^l+\tfrac{\Lambda_1^l}{2}-\lambda_2\right)
\tfrac{\Lambda_2^l+2R_2^l-2\lambda_2}{\Lambda_2^l}
\right]e^{-\lambda_2 y}.
\label{soln:tilde-p1-m}
\end{align}
The bounded solution for $\widetilde p_1^\pm$ is determined by the free parameters $\widetilde C_1(t)$ and $\widetilde C_2(t)$. The $\tilde p_0^\pm$-equations read
\begin{align*}
\partial_y\widetilde p_0^+= R_1^l\widetilde p_1^+,
\qquad
\partial_y\widetilde p_0^-
= -R_1^l\widetilde p_1^-.
\end{align*}
Since $\widetilde p_1^\pm$ decay exponentially as $y\to\infty$, the right-hand sides are integrable on $(0,\infty)$. Hence there exist $\widetilde p_0^\pm(\infty)$ such that
\begin{align}
\widetilde p_0^+(t,y)
= &\, \widetilde p_0^+(\infty)
- \tfrac{R_1^l\Lambda_1^l}{2\lambda_1}\widetilde C_1
e^{-\lambda_1 y} \no
\\
&- \tfrac{R_1^lR_2^l\widetilde C_2}{\lambda_2 ((R_1^l)^2+R_1^l\Lambda_1^l-\lambda_2^2)}
[R_1^l+\tfrac{\Lambda_1^l}{2}+\lambda_2
+\tfrac{\Lambda_1^l}{2}
\tfrac{\Lambda_2^l+2R_2^l-2\lambda_2}{\Lambda_2^l}]e^{-\lambda_2 y},
\label{soln:tilde-p-1-p}
\\
\widetilde p_0^-(t,y)
= & \,
\widetilde p_0^-(\infty)
+ \tfrac{R_1^l}{\lambda_1}
(R_1^l+\tfrac{\Lambda_1^l}{2}-\lambda_1)
\widetilde C_1
e^{-\lambda_1 y} \no
\\
& + \tfrac{R_1^lR_2^l\widetilde C_2}{\lambda_2 ((R_1^l)^2+R_1^l\Lambda_1^l-\lambda_2^2)}
[\tfrac{\Lambda_1^l}{2}
+(R_1^l+\tfrac{\Lambda_1^l}{2}-\lambda_2)
\tfrac{\Lambda_2^l+2R_2^l-2\lambda_2}{\Lambda_2^l}]e^{-\lambda_2 y}.
\label{soln:tilde-p-1-m}
\end{align}
In the case where $\lambda_1 = \lambda_2 = \lambda$, the solutions for $\tilde p_1^\pm(t,y), \tilde p_0^\pm(t,y)$ are as follows:
\begin{align*}
\widetilde p_1^+(t,y)
=&\, \tfrac{\Lambda_1^l}{2}\widetilde C_1(t)e^{-\lambda y}
+\tfrac{R_2^l\widetilde C_2(t)}{2\lambda}
[R_1^l+\tfrac{\Lambda_1^l}{2}+\lambda
+\tfrac{\Lambda_1^l}{2}\tfrac{\Lambda_2^l+2R_2^l-2\lambda}{\Lambda_2^l}] y e^{-\lambda y}
\\
& \quad \,
+\tfrac{R_2^l\widetilde C_2(t)}{4\lambda^2}
[R_1^l+\tfrac{\Lambda_1^l}{2}-\lambda
+\tfrac{\Lambda_1^l}{2}\tfrac{\Lambda_2^l+2R_2^l-2\lambda}{\Lambda_2^l}] e^{-\lambda y},
\\ 
\widetilde p_1^-(t,y)
=&\, (R_1^l+\tfrac{\Lambda_1^l}{2}-\lambda)
\widetilde C_1(t)e^{-\lambda y}
+\tfrac{R_2^l\widetilde C_2(t)}{2\lambda}
[\tfrac{\Lambda_1^l}{2}
+
\tfrac{(R_1^l+\tfrac{\Lambda_1^l}{2}-\lambda)(\Lambda_2^l+2R_2^l-2\lambda)}{\Lambda_2^l}] y e^{-\lambda y}
\\
& \quad \,
+ \tfrac{R_2^l\widetilde C_2(t)}{4\lambda^2}
[\tfrac{\Lambda_1^l}{2}
+ \tfrac{(R_1^l+\tfrac{\Lambda_1^l}{2}+\lambda)(\Lambda_2^l+2R_2^l-2\lambda)}{\Lambda_2^l}] e^{-\lambda y}, 
\end{align*}
and
\begin{align*}
\widetilde p_0^+(t,y)
=&\, \widetilde p_0^+(\infty)
- \tfrac{R_1^l\Lambda_1^l}{2\lambda}\widetilde C_1 e^{-\lambda y}
\\
& - \tfrac{R_1^lR_2^l\widetilde C_2}{2\lambda}
[
R_1^l+\tfrac{\Lambda_1^l}{2}+\lambda
+\tfrac{\Lambda_1^l}{2}
\tfrac{\Lambda_2^l+2R_2^l-2\lambda}{\Lambda_2^l}]
e^{-\lambda y}
(\tfrac{y}{\lambda}+\tfrac{1}{\lambda^2})
\\
& - \tfrac{R_1^lR_2^l\widetilde C_2}{4\lambda^3}
[
R_1^l+\tfrac{\Lambda_1^l}{2}-\lambda
+\tfrac{\Lambda_1^l}{2}
\tfrac{\Lambda_2^l+2R_2^l-2\lambda}{\Lambda_2^l}]
e^{-\lambda y},
\\
\widetilde p_0^-(t,y)
=&\, \widetilde p_0^-(\infty)
+\tfrac{R_1^l}{\lambda}
(R_1^l+\tfrac{\Lambda_1^l}{2}-\lambda)
\widetilde C_1 e^{-\lambda y}
\\
&+ \tfrac{R_1^lR_2^l\widetilde C_2}{2\lambda}[
\tfrac{\Lambda_1^l}{2}
+(R_1^l+\tfrac{\Lambda_1^l}{2}-\lambda)
\tfrac{\Lambda_2^l+2R_2^l-2\lambda}{\Lambda_2^l}]
e^{-\lambda y}(\tfrac{y}{\lambda}+\tfrac{1}{\lambda^2})
\\
&+ \tfrac{R_1^lR_2^l\widetilde C_2}{4\lambda^3}[
\tfrac{\Lambda_1^l}{2}
+(R_1^l+\tfrac{\Lambda_1^l}{2}+\lambda)
\tfrac{\Lambda_2^l+2R_2^l-2\lambda}{\Lambda_2^l}]
e^{-\lambda y}.
\end{align*}
In both cases, it is clear that the left-domain solution operator are determined by the asymptotic data $(\widetilde p_0^+(\infty),\widetilde p_0^-(\infty))^\top$.

At $y=0$, if $\lambda_1 \neq \lambda_2$, then we have
\begin{align} \label{interf:tilde-p}
&\widetilde p_0^+(t,0)
=  \widetilde p_0^+(\infty)
- \tfrac{R_1^l\Lambda_1^l}{2\lambda_1}\widetilde C_1  \no
\\
& \hspace{2cm}
- \tfrac{R_1^lR_2^l\widetilde C_2}{\lambda_2 ((R_1^l)^2+R_1^l\Lambda_1^l-\lambda_2^2)}
[R_1^l+\tfrac{\Lambda_1^l}{2}+\lambda_2
+\tfrac{\Lambda_1^l}{2}
\tfrac{\Lambda_2^l+2R_2^l-2\lambda_2}{\Lambda_2^l}], \no
\\
&\widetilde p_0^-(t,0)
=
\widetilde p_0^-(\infty)
+ \tfrac{R_1^l}{\lambda_1}
(R_1^l+\tfrac{\Lambda_1^l}{2}-\lambda_1)
\widetilde C_1 \no
\\
& \hspace{2cm} + \tfrac{R_1^lR_2^l\widetilde C_2}{\lambda_2 ((R_1^l)^2+R_1^l\Lambda_1^l-\lambda_2^2)}
[\tfrac{\Lambda_1^l}{2}
+(R_1^l+\tfrac{\Lambda_1^l}{2}-\lambda_2)
\tfrac{\Lambda_2^l+2R_2^l-2\lambda_2}{\Lambda_2^l}],   
\\
&\widetilde p_1^+(t,0)
=\tfrac{\Lambda_1^l}{2}
\widetilde C_1(t) \no
+\tfrac{R_2^l\widetilde C_2(t)}
{(R_1^l)^2+ R_1^l\Lambda_1^l -\lambda_2^2}
\left[
R_1^l+\tfrac{\Lambda_1^l}{2}+\lambda_2
+\tfrac{\Lambda_1^l}{2}
\tfrac{\Lambda_2^l+2R_2^l-2\lambda_2}{\Lambda_2^l}
\right], \no
\\
&\widetilde p_1^-(t,0)
=\left(R_1^l+\tfrac{\Lambda_1^l}{2} - \lambda_1\right)
\widetilde C_1(t) \no
\\
& \hspace{2cm}
+\tfrac{R_2^l\widetilde C_2(t)}
{(R_1^l)^2+ R_1^l\Lambda_1^l-\lambda_2^2}
\left[
\tfrac{\Lambda_1^l}{2}
+\left(R_1^l+\tfrac{\Lambda_1^l}{2}-\lambda_2\right)
\tfrac{\Lambda_2^l+2R_2^l-2\lambda_2}{\Lambda_2^l}
\right], \no
\\
&\widetilde p_2^+(t,0)=\widetilde C_2(t),
\qquad
\widetilde p_2^-(t,0)
=\frac{\Lambda_2^l+2R_2^l-2\lambda_2}{\Lambda_2^l}
\widetilde C_2(t).    \no 
\end{align}

\medskip
\noindent\textbf{Right domain $\Omega_r$} In the right domain, equation (3.8) becomes
\begin{align}\label{sys:mn-leading order}
  \Omega_{r}: \quad  
\begin{cases}
    \p_z  \widehat p_0^+ -  R_0^r  \widehat p_0^+ = 0,
\\[2mm]
     \p_z  \widehat p_0^- +  R_0^r  \widehat p_0^- = 0,
\\[2mm]
     \p_z  \widehat p_1^+ + R_{0}^r  \widehat p_{0}^+ - R_1^r  \widehat p_1^+ = 0,
\\[2mm]
    \p_z  \widehat p_1^- -  R_{0}^r  \widehat p_{0}^- +  R_{1}^r  \widehat p_{1}^- = 0,
\\[2mm]
    \p_z  \widehat p_2^+ + R_{1}^r  \widehat p_{1}^+ = 0,
\\[2mm]
    \p_z  \widehat p_2^- -  R_{1}^r  \widehat p_{1}^- = 0.
\end{cases}
\end{align}
Given the structure of system \eqref{sys:mn-leading order} and the requirement that the solutions $\widehat{p_{i}}^\pm$ $(i=0,1,2)$ remain bounded as $z \to \infty$, we solve the first two equations of the system yields
\begin{align}\label{solu:right-0}
\widehat p_0^+(t,z)= 0, 
\qquad 
  \widehat p_0^-(t,z)= \widehat C_0(t) e^{-R_0^r z},
\end{align}
where $\widehat C_0(t)$ is a free parameter. 

Similar as for the left domain, we separate the two cases of $R^r_0 = R^r_1$ and $R^r_0 \neq R^r_1$. If $R^r_0 \neq R^r_1$, then by substituting the solution into the equation for $\widehat p_1^\pm$, we obtain
\begin{align}
 \widehat p_1^+ (t,z) = 0, 
\qquad 
  \widehat p_1^- (t,z) = \tfrac{R_0^r\widehat C_0}{R_1^r-R_0^r}e^{-R_0^r z}
+\widehat C_1 e^{-R_1^r z},
\label{soln:hat-p-1}
\end{align}
From the equation for $\widehat p_2^\pm$, since $\widehat p_1^\pm$ decay exponentially as $z\to\infty$, there exist $\widehat p_2^\pm(\infty)$ such that
\begin{align}
    \widehat p_2^+(t,z)=\widehat p_2^+(\infty),
\quad
  \widehat p_2^-(t,z)=\widehat p_2^-(\infty) -\tfrac{R_1^r \widehat C_0}{R_1^r - R_0^r}e^{-R_0^rz} -\widehat C_1 e^{-R_1^r z},
\label{soln:hat-p-2}
\end{align}
where $R_1^r\neq R_0^r$. When $R_1^r = R_0^r = R$, we have
\begin{align*}
    \widehat p_1^+(t,z) &= 0, \quad \widehat p_1^-(t,z)= \bigl(R\widehat C_0 z+\widehat C_1\bigr)e^{-Rz},
\\
\widehat p_2^+(t,z) &= \widehat p_2^+(\infty),
\quad
\widehat p_2^-(t,z) = \widehat p_2^-(\infty) +\bigl[-R\widehat C_0 z-\widehat C_0-\widehat C_1\bigr]e^{-Rz}.
\end{align*}

At $z=0$, if $R^r_0 \neq R^r_1$, then we have
\begin{align} \label{interf:hat-p}
&\widehat p_0^+(t,0)= 0, 
\qquad
\widehat p_0^-(t,0)= \widehat C_0(t),  \no
\\
& \widehat p_2^+(t,0)=\widehat p_2^+(\infty),
\quad
  \widehat p_2^-(t,0)=\widehat p_2^-(\infty) -\tfrac{R_1^r \widehat C_0}{R_1^r - R_0^r} -\widehat C_1,
\\
& \widehat p_2^+(t,0)=\widehat p_2^+(\infty),
\quad
  \widehat p_2^-(t,0)=\widehat p_2^-(\infty) -\tfrac{R_1^r \widehat C_0}{R_1^r - R_0^r} -\widehat C_1. \no
\end{align}
Using~\eqref{interf:hat-p}, \eqref{interf:tilde-p} and the continuity equation, we eliminate $\widetilde C_1, \widetilde C_2, \widehat C_0, \widehat C_1$ and obtain the interface relation between $\widetilde p_0^\pm(\infty)$ and $\widehat p_2^\pm(\infty)$. Recall that these density functions are those with $v \to -v$. So we change the sign of $v$ back and obtain the matrix $A_2$ given in (3.23).

\bigskip
\bigskip
\section{Algorithm for continuous $m$}

\quad
 
\begin{algorithm}[H] 
\caption{Particle motion with continuous $m$ and discrete $v=\pm v_0$}
\label{alg: continuous mv0}
\begin{algorithmic}[1]
\Require 
\State Computational domain $[X_l,X_r]$; terminal time $T_t$; external environment $m_c(x)$; time step $\Delta t$.
\State For the $j$-th particles:  
 initial position $x_{(j)}$; initial internal state $m_{(j)}$; initial velocity $v_{(j)}$; 
\Ensure  
\Repeat
\For{each $j$}
\State Determine $\lambda_{(j)}=\Lambda(m_{(j)}-m_c(x_{(j)}))$
    \State Generate one random number $r$ uniformly distributed in $[0, 1]$:  
    \If{$r \leq \frac{1}{2} (1 - \exp(-\lambda_{(j)} \Delta t))$}
        \State $v_{(j)} \gets -v_{(j)}$  
    \Else
        \State $v_{(j)} \gets v_{(j)}$  
    \EndIf
    \State $x_{(j)} \gets x_{(j)} + v_{(j)}\Delta t$;  
    \If{$x_n > X_r$}
        \State $x_{(j)} \gets x_{(j)}-|X_r-X_l|$
    \EndIf
    \If{$x_{(j)} < X_l$}
        \State $x_{(j)} \gets x_{(j)}+|X_r-X_l|$
    \EndIf
    \State $m_{(j)} \gets m_{(j)} +  \Delta t\cdot f(m_{(j)}-m_c(x_{(j)}))$; 
    \State $\lambda_{(j)} \gets \Lambda(m_{(j)},m_c(x_{(j)}))$;  
\EndFor
\State $T \gets T + \Delta t$  
\Until{$T = T_t$}
\end{algorithmic} 
\end{algorithm}

\newpage
\section{Algorithm for the multi-species model}

\quad

\begin{algorithm}[htbp]
\caption{Particle motion with discrete $m$ and discrete $v=\pm v_0$}
\label{alg: discrete mv}
\begin{algorithmic}[1]
\Require
\State Computational domain $[X_l,X_r]$; terminal time $T_t$; subdomains $\Omega_l$, $\Omega_r$; time step $\Delta t$; grid spacing $\Delta m = (M_r-M_l)/n$.
\State For the $j$-th particle: initial position $x_{(j)}$, initial species label $i_{(j)}\in\{0,\dots,n\}$, initial velocity $v_{(j)}\in\{-v_0,+v_0\}$.
\Ensure
\Repeat
\For{each particle $j$}
    \State Update $v_{(j)}$ and $x_{(j)}$ as in Algorithm \ref{alg: continuous mv0}
    
    \State Set $m_c = M_l$ if $x_{(j)}\in\Omega_l$, else $m_c = M_r$
    \State $m_{\text{half}} \gets m_{i_{(j)}} - \Delta m/2$ (if $\Omega_l$) else $m_{i_{(j)}} + \Delta m/2$
    \State $f \gets \frac{k_R}{\varepsilon}\bigl(1 - a(m_{\text{half}}-m_c)/a_0\bigr)$
    \If{$(x_{(j)}\in\Omega_l \text{ and } f < 0)$ or $(x_{(j)}\in\Omega_r \text{ and } f > 0)$}
        \State $R \gets f / \Delta m$
    \Else
        \State $R \gets 0$
    \EndIf
    \State Generate $p\sim\mathcal{U}[0,1]$
    \If{$p \le 1-\exp(-|R|\,\Delta t)$}
        \If{$R < 0$}
            \State $i_{(j)} \gets i_{(j)}-1$
        \Else
            \State $i_{(j)} \gets i_{(j)}+1$
        \EndIf
    \EndIf
    \State $i_{(j)} \gets \max(0,\;\min(n,\;i_{(j)}))$
    \State Set $m_{(j)} \gets m_{i_{(j)}}$ 
    \State $\lambda_j \gets \Lambda\big(m_{(j)} - m_c(x_{(j)})\big)$ 
\EndFor
\State $T \gets T + \Delta t$
\Until{$T = T_t$}
\end{algorithmic}
\end{algorithm}